\documentclass[paper=A4, fontsize=14pt, headings=normal, abstract=on, headsepline, DIV=calc]{scrartcl}

\usepackage[utf8]{inputenc} 
\usepackage[T1]{fontenc}
\usepackage{lmodern}
\usepackage{microtype}

\usepackage{amsxtra,amsthm}
\usepackage{amssymb,latexsym}
\usepackage{upref}
\usepackage{exscale}

\usepackage{tikz-cd}
\usepackage{scalerel}
\usepackage{mathtools}

\usepackage[english]{babel}

\usepackage[pdftitle={[Accepted manuscript] Arason's filtration of the Witt group of dyadic valued fields},
pdfauthor={Joachim Verstraete}]{hyperref}

\usepackage[abbrev,backrefs]{amsrefs}

\usepackage{scrlayer-scrpage}
\lohead*{\small This is the accepted manuscript. © 2019. This manuscript version is made available under the CC-BY-NC-ND 4.0 license \url{http://creativecommons.org/licenses/by-nc-nd/4.0/}.  \\
 DOI of the published version  : DOI \href{https://doi.org/10.1016/j.jalgebra.2018.10.026}{10.1016/j.jalgebra.2018.10.026}.}
\KOMAoptions{DIV=last}

\numberwithin{equation}{section}

\newtheorem{theorem}{Theorem}[section]
\newtheorem{proposition}[theorem]{Proposition}	
\newtheorem{corollary}[theorem]{Corollary}	
\newtheorem{lemma}[theorem] {Lemma}		
\newtheorem{definition}[theorem] {Definition}

\theoremstyle{definition}
\newtheorem{example}[theorem] {Example}

\theoremstyle{remark}
\newtheorem{remark}[theorem] {Remark}

\newcommand*{\defn}[2][]{\begin{definition}[#1]\emph{ #2 }\end{definition} }
\newcommand*{\lem}[2][]{\begin{lemma}[#1] #2 \end{lemma} }
\newcommand*{\prop}[2][]{\begin{proposition}[#1] #2 \end{proposition} }
\newcommand*{\cor}[2][]{\begin{corollary}[#1] #2 \end{corollary}}
\newcommand*{\thm}[2][]{\begin{theorem}[#1] #2 \end{theorem}}

\newcommand*{\rem}[2][]{\begin{remark}[#1] #2 \end{remark} }

\newcommand*{\pr}[1]{\begin{proof} #1 \end{proof}}

\newcommand*{\listeN}[1]{\begin{enumerate} #1 \end{enumerate}}

\let\ensembleNombre\mathbb
\newcommand*\N{\ensuremath{\ensembleNombre{N}}}
\newcommand*\Z{\ensuremath{\ensembleNombre{Z}}}
\newcommand*\Q{\ensuremath{\ensembleNombre{Q}}}

\newcommand*{\set}[1]{\left\{ #1 \right\} }
\newcommand*{\tq}{\,\vert \,}

\newcommand*\iso{\ensuremath{\xrightarrow{\,\,\,\,\sim\,\,\,\,}}}

\DeclareMathOperator{\gr}{gr}
\DeclareMathOperator{\sev}{span}
\DeclareMathOperator{\aut}{Aut}

\DeclareMathOperator{\charact}{char}
\DeclareMathOperator{\charac}{char}
\DeclareMathOperator*{\uperp}{\scalerel*{\perp}{\sum}}

\newcommand*{\mP}{ {\mathcal{P}} }
\newcommand*{\mPp}{ {\mathcal{P}_{p}} }
\newcommand*{\mPm}{ {\mathcal{P}_{m}} }

\begin{document}

  \title{Arason's filtration of the Witt group of dyadic valued fields.}
  \author{Joachim \textsc{Verstraete}\thanks{Joachim Verstraete is a Research Fellow of the Fonds de la Recherche Scientifique - FNRS.}}

  \date{\today}
  
  \maketitle
    
  \begin{abstract}
  Generalizing a theorem of Springer, we construct an extended Arason filtration by subgroups for the quadratic Witt group of a general valued field, relating these subgroups with Witt-like groups of the residue field, in arbitrary characteristic. These Witt-like groups involve totally singular quadratic forms. In the case of a discretely valued field, we recover the original Arason filtration.
  

  \medskip\noindent
  \textsc{Classification (MSC 2010):} 11E81

  \end{abstract}
  
  \tableofcontents
  
%%%%%%%%%%%%%%%%%%%%%%%%%%%%%%%%%%%%%%%%%
\vspace{1em}

A well-known result of Springer states that the Witt group of quadratic forms \(W_q(F)\) of a complete discretely valued field \(F\)  is isomorphic to a direct sum of two copies of the Witt group of quadratic forms \(W_q(\overline{F})\) of the residue field \
\(\overline{F}\), provided that the residue characteristic is different from \(2\) (see \cite{MR0070664}). Subsequently, that result has been generalized to Henselian valued fields \(F\) of residue characteristic different from \(2\), see for example the paper of 
Tietze~\cite{MR0366809}. Tietze also considered the case when \(\charact \overline{F} = 2\), and he obtained a similar isomorphism from a subgroup $U(F) \subset W_q(F)$ onto a direct sum of copies of \(W_q(\overline{F})\). In their paper~\cite{MR1385276} 
in 1996, Aravire and Jacob performed an exhaustive analysis of the Witt group of a dyadic maximally complete field with perfect residue field and showed that the description of the Witt group in that case is extremely delicate. More recently, Arason proposed, 
in his note \cite{RH-05-2016}, another way to describe completely the Witt group of a discretely valued field by a filtration by subgroups, covering the case of residue characteristic 2. Our purpose in this paper is to obtain results similar to  Arason's, but for 
general valued fields, by using techniques of \cite{MR2836073} involving graded structures arising from the valuation. Our approach has the main advantage to employ more intrinsic methods than the computational ones of \cite{RH-05-2016}.  In particular, we 
develop special Witt-like groups involving totally singular quadratic forms.

More precisely, our main results are the following (expressed in a similar way as \cite{RH-05-2016}):

Let $F$ be a field with a (not necessarily surjective) valuation $v \colon F \to \Gamma \cup \set{\infty}$, where $\Gamma$ is a totally ordered abelian group. Without loss of generality, we assume that $\Gamma$ is divisible, since we may substitute for $
\Gamma$ its divisible hull. Denote by $\Gamma_F$ the value group of $(F,v)$, i.e. \(\Gamma_F = v(F^{\times}) \subset \Gamma\), and let \(\overline{F}\) denote the residue field of \((F,v)\). Let $W_q(F)$ (resp. $W(F)$) be the Witt group of (nonsingular) 
quadratic forms (resp. the Witt group of (nondegenerate) symmetric bilinear forms) over $F$. Recall that $W_q(F)$ can be identified with $W(F)$ when $2$ is invertible in $F$ (to a symmetric bilinear form $b$, we associate the quadratic form given by $x 
\mapsto b(x,x))$. \\

1) It is well-known (see for example \cite{MR2836073}*{Proposition 8 and Corollary 11}) (but a proof of this fact will also be provided in the paper) that there always exists a subgroup $W_q(F)_{tame} \subset W_q(F)$ together with a group epimorphism 
\begin{equation}\label{springer}
  W_q(F)_{tame} \to \bigoplus_{\Gamma_F / 2 \Gamma_F} W_q(\overline{F}),
\end{equation}
which is an isomorphism when $F$ is Henselian. That epimorphism depends on a choice of uniformizing parameters. When $\charact(\overline{F}) \neq 2$, we always have $W_q(F)_{tame} = W_q(F)$  and, when moreover $F$ is Henselian, the 
isomorphism~\eqref{springer} is Springer's theorem. When $\charact(\overline{F}) = 2$, it turns out that $W_q(F)_{tame}$ equals the previously mentioned Tietze subgroup $U(F)$, as shown in \cite{MR2836073}*{Corollary 17}. \\

Let us now give some informal motivation for the next, but principal, result. If \(F\) is a dyadic (i.e., \(\charact \overline{F} = 2\)) Henselian valued field, the preceding isomorphism gives a description of \(W_q(F)_{tame}\). Therefore, it remains to describe the 
quotient \({W_q(F)}/{W_q(F)_{tame}}\) (this quotient is called \(W_q(F)_{wild}\) by Arason in \cite{RH-05-2016} when \(F\) is a discretely valued field). It is the description of the structure of the wild part which motivates the following result. (Note that the result 
does not require \(F\) to be Henselian).\\

2) (Main result) If $\charact(\overline{F}) = 2$, there exists an ascending filtration
\[
    \begin{aligned}
        &\left(W_q\left(F\right)_{\varepsilon}\right)_{\varepsilon \in E}, &E = \set{\varepsilon \in \frac 1 2 \Gamma_F \tq 0 \leq \varepsilon \leq v(2)}
    \end{aligned}
\]
 of $W_q(F)$ by subgroups $W_q(F)_{\varepsilon} \subset W_q(F)$.
This filtration is called the \emph{Arason filtration of $W_q(F)$} and satisfies
\[
    \begin{aligned}
        &W_q(F)_0 = W_q(F)_{tame} &\text{ and } &\bigcup_{\varepsilon \in E}W_q(F)_{\varepsilon} = W_q(F).
    \end{aligned}
\]
Moreover, for each $\varepsilon \in \frac 1 2 \Gamma_F$ such that $0 < \varepsilon \leq v(2)$,
there is a group isomorphism
\begin{equation}\label{filtration_isos}
  \frac{W_q(F)_{\varepsilon}} {W_q(F)_{< \varepsilon}} \iso
  \begin{cases}
    \bigoplus_{\frac {\left| \Gamma_F / 2 \Gamma_F \right|} 2} W_{ssq}(\overline{F}) & \text{if } \varepsilon \notin \Gamma_F \\
    \bigoplus_{\Gamma_F / 2 \Gamma_F} W_{sq}(\overline{F}) & \text{if } \varepsilon \in \Gamma_F  \text{ and } \varepsilon \neq v(2) \\
    \bigoplus_{\Gamma_F / 2 \Gamma_F} W(\overline{F}) & \text{if } \varepsilon = v(2)
  \end{cases}
\end{equation}
where \[W_q(F)_{< \varepsilon } = \bigcup_{\substack{\gamma \in E, \\ 0 \leq \gamma < \varepsilon}} W_q(F)_{\gamma},\] and $W_{sq}(\overline{F})$ along with $W_{ssq}(\overline{F})$ are two Witt-like groups of $\overline{F}$, with $W_{ssq}(\overline{F})$ 
nontrivial (see point $3$ below). These isomorphisms depend on a choice of uniformizing parameters. Note that if $\charact(F) = 2$, then $v(2) = \infty$, hence $E = \set{\varepsilon \in \frac 1 2 \Gamma_F \tq \varepsilon \geq 0} $. Note also that if $\charact(F) 
\neq 2$, then we have \(v(2) \in \Gamma_F\) and $W_q(F)_{v(2)} = W_q(F)$ (since the filtration is increasing and covers \(W_q(F)\)).

In particular, if $\Gamma_F$ is a well-ordered group (i.e. $\Gamma_F \simeq \Z$ as ordered groups), then for each $\varphi \in W_q(F)$ there is a minimal $\varepsilon \geq 0$ such that $\varphi \in W_q(F)_{\varepsilon}$. Under the same hypothesis, we also 
have that $W_q(F)_{< \varepsilon} = W_q(F)_{\gamma}$ for some $\gamma \in E$. \\

Moreover, for \(a,b \in F\), denote by \([a,b]\) the two-dimensional quadratic form \(q \colon F \times F \to F\) given by \(q(x_1,x_2) = ax_1^2 + x_1x_2 + bx_2^2\).

Corollary~\ref{W_q__asGenerated} states that, for each \(\varepsilon \in \frac 1 2 \Gamma_F\) such that \( 0 \leq \varepsilon < v(2) \), \(W_q(F)_{\varepsilon}\) is the subgroup of \(W_q(F)\) generated by the Witt classes of two-dimensional quadratic forms \([a,b]
\) with \(a, b\in F\) satisfying \(v(a) + v(b) \geq - 2 \varepsilon\). The condition \(\varepsilon < v(2)\) ensures that \([a,b]\) is nonsingular. Note also that if \(F\) is Henselian and \(v(a) + v(b) > 0\), the proof of Lemma~\ref{calculations_discrVal} \((e)\) shows that the 
Witt class of \([a,b]\) is \(0\).

\begin{example}
\label{ex_intro_1}
	Let \(k\) be a field of characteristic \(2\). Let \(F = k((t))\) be the field of formal Laurent series in an indeterminate \(t\) with coefficients in \(k\). 
	Equip \(F\) with the usual \(t\)-adic valuation \(v \colon F \to \Z \cup \set{\infty}\), so that \(v(t) = 1\) and \(v(\alpha) =0\) for all \(\alpha \in k \subset F\).
	Let also \(q_W\) be the Witt class of a quadratic form \(q\) in \(W_q(F)\). It can be shown (by using Lemma~\ref{calculations_discrVal} \((d)\)) that the induced map \([ \phantom{x} \, , \, \phantom{x}]_W \colon F \times F \to W_q(F)\) is bi-additive since \
\(\charact F = 2\). 
	
	For this field \(F\), we have an infinite ascending filtration of \(W_q(F)\) indexed by \(\frac 1 2 \N\) : 
	\[
	 W_q(F)_{tame} = W_q(F)_0 \subset W_q(F)_{1/2} \subset W_q(F)_{1} \subset W_q(F)_{3/2} \subset ... \, ,
	\]
	and the form \([a,b]_W\) (with \(a,b \in F\)) belongs to \(W_q(F)_{\varepsilon} \) if (but not only if) \(v(a) + v(b) \geq -2 \varepsilon\), for all \(\varepsilon \in \frac 1 2 \N\).
	 
	Consider for example the quadratic form \([1 + t, t^{-1} + t]\). By bi-additivity of \([ \phantom{x} \, , \, \phantom{x}]_W\), we have
	\[
	      [1 + t, t^{-1} + t]_W = [1, t^{-1}]_W + [t, t^{-1}]_W + [1,t]_W + [t,t]_W.
	\]
	Since \(v(1) + v(t) = 1 > 0\) and \(F\) is Henselian, it follows that \([1,t]_W = 0\). Similarly, \([t,t]_W = 0\).
	We have \([t, t^{-1}] _W \in W_q(F)_{tame}\) (since \(v(t) + v(t^{-1}) \geq 0\)), and \([1, t^{-1}]_W \in W_q(F)_{1/2}\) (since \(v(1) + v(t^{-1}) \geq -1 = -2 \cdot \frac 1 2\)). But it turns out that \([1, t^{-1}]_W \not \in W_q(F)_{tame}\). (It is an application of 
Corollary~\ref{Arason_2}.)
	Therefore
	\[
	  [1 + t, t^{-1} + t]_W = [1, t^{-1}]_W + [t, t^{-1}]_W  \in W_q(F)_{1/2} \setminus W_q(F)_{0}.
	\]
	
	Suppose now that \(\alpha \in k^{\times} \setminus {k^{\times}}^{2}\). Then Corollary~\ref{Arason_2} shows as previously that \([1, \alpha t^{-2}]_W \in W_q(F)_{1} \setminus W_q(F)_{1/2} \) ;
	but note that \([1, t^{-2}] \simeq [1, t^{-1}]\) (by the change of basis given by \(e_1' = e_1, e_2' = e_2 + t^{-1} e_1\), where \(e_1, e_2\) is the standard basis of \([1, t^{-2}]\)), so that \([1, t^{-2}]_W = [1, t^{-1}]_W  \in W_q(F)_{1/2}\). \qed
\end{example}
\begin{example}
	This example uses the same notations as Example~\ref{ex_intro_1}. The residue field \(\overline{F}\) of \(F = k((t))\) for the \(t\)-adic valuation \(v\) is canonically isomorphic to \(k\). \\
	
	Consider the nonsingular quadratic form \([t, t^{-1}]\) over \(F\). By Lemma~\ref{calculations_discrVal} \((c)\), that form is isometric to the form \(q_1 \coloneqq t[1,1]\), that is the quadratic form given for all \((x_1, x_2) \in F \times F\) by
	\[
	      q_1(x_1,x_2) = t(x_1^2 + x_1x_2 + x_2^2).
	\]
	Its nondegenerate associate polar bilinear form \(b_{q_1}\) is given for all \((x_1, x_2), (y_1,y_2) \in F \times F\) by
	\[
	      b_{q_1}((x_1,x_2),(y_1,y_2)) = t(x_1y_2 + x_2 y_1).
	\]
	
	Under the group isomorphism \(W_q(F)_{tame} \to W_q(k) \oplus W_q(k)\), the Witt class \(q_{1,W}\) is sent to \((0, q'_{1,W})\), where \(q'_1\) is the nonsingular quadratic form over \(k\) given for all \((x_1, x_2) \in k \times k\) 
	\[
	      q'_1(x_1,x_2) = x_1^2 + x_1x_2 + x_2^2.
	\]
	Its nondegenerate associate polar bilinear form \(b_{q'_1}\) is given for all \((x_1, x_2), (y_1,y_2) \in k \times k\) by
	\[
	      b_{q'_1}((x_1,x_2),(y_1,y_2)) = x_1y_2 + x_2 y_1.
	\]
	Similarly, under the same map, the Witt class \([1,1]_W\) is sent to \((q'_{1,W}, 0)\). (See Theorem~\ref{filtration} and Proposition~\ref{decompsitionCan_1} (or Corollary~\ref{iso_gr_sq_1}) for the exact description of the map.)\\
	
	Consider now the nonsingular quadratic form \([1, \alpha t^{-2}]\) for some arbitrary \(\alpha \in k^{\times} \setminus {k^{\times}}^2\). It is isometric to the quadratic form \(q_2\) given for all \((x_1, x_2) \in F \times F\) 
	\[
	      q_2(x_1,x_2) = x_1^2 + t x_1x_2 + \alpha x_2^2.
	\]
	Its nondegenerate associate polar bilinear form \(b_{q_2}\) is given for all \((x_1, x_2), (y_1,y_2) \in F \times F\) by
	\[
	      b_{q_2}((x_1,x_2),(y_1,y_2)) =t(x_1y_2 + x_2 y_1).
	\]
	Under the group isomorphism \(W_q(F)_{1} \to W_{sq}(k) \oplus W_{sq}(k)\), the Witt class \(q_{2,W}\) is sent to \((\Phi_W, 0)\), where \(\Phi_W\) is the Witt class in \(W_{sq}(k)\) of a \emph{symplectic quadratic space} \(\Phi = (q'_2, b'_{q_2})\) (see 
Definition~\ref{def_ssq} for a precise definition). That space is made up, partly, of
	 a quadratic form \(q'_2\) over \(k\), induced by \(q_2\), which is given for all \((x_1, x_2) \in k \times k\) by
	\[
	      q'_2(x_1,x_2) = x_1^2  + \alpha x_2^2.
	\]
	Therefore, \(q'_2\) is a \emph{totally singular} quadratic form.
	But the bilinear form \(b_{q_2}\) yields, after scaling by \(t^{-1}\), the \emph{nondegenerate} alternating bilinear form \(b'_{q_2}\) over \(k\) given for all \((x_1, x_2), (y_1,y_2) \in k \times k\) by
	\[
	      b'_{q_2}((x_1,x_2),(y_1,y_2)) = x_1y_2 + x_2 y_1.
	\]
	Similarly, under the same map, the Witt class of the form \(t[1, \alpha t^{-2}]\) is sent to \((0, \Phi_{W})\). (See Theorem~\ref{filtration} and Proposition~\ref{decompsitionCan_1} (or Corollary~\ref{iso_gr_sq_1}) for the exact description of the map.)\\
	
	Note also that, by Lemma~\ref{calculations_discrVal} \((c)\), the quadratic form \([1, t^{-1}]\), which satisfies \([1, t^{-1}]_W \in W_q(F)_{1/2} \setminus W_q(F)_{0}\), is such that \(t[1, t^{-1}] \simeq [1, t^{-1}]\). \qed
\end{example}
\begin{example}
	Let \(F = \Q_2\) be the field of dyadic numbers equipped with the usual \(2\)-adic valuation, so that \(v(2) = 1\). Its residue field is \(\overline{F} = \mathbb{F}_2\), the finite field with two elements. In this case, we have a \(3\)-step filtration 
	\[
	W_q(F)_{tame} = W_q(F)_0 \subset W_q(F)_{1/2} \subset W_q(F)_{1} = W_q(F).
	\]
	
	Consider the one-dimensional nonsingular quadratic form \(q\) over \(F\) given for all \(x \in F\) by 
	\[
	     q(x) = x^2,
	\]
	which can be identified (since \(\charact F \neq 2 \)) with the nondegenerate symmetric bilinear form \(b\) given for all \(x,y \in F\) by
	\[
	    b(x,y) = xy.
	\]
	Under the group isomorphism \(\partial^{1} \colon W_q(F)_{1} \to W(\overline{F}) \oplus W(\overline{F})\), the Witt class \({q}_W\) is sent to \(({b'}_W,0)\), where \(b'\) is the nondegenerate one-dimensional symmetric bilinear form over \(\mathbb{F}_2\) 
given for all \(x, y \in \mathbb{F}_2\) by
	\[
	      b'(x,y) =xy.
	\]
	Its associate totally singular quadratic form \(q_{ b'}\) over \(\mathbb{F}_2\) is given for all \(x \in \mathbb{F}_2\) by 
	\[
	 q_{ b'}(x) = x^2.
	\]
	Similarly, under the same map, the Witt class of the form \(2q\) is sent to \((0, b'_W)\). (See Theorem~\ref{filtration} and Proposition~\ref{decompsitionCan_1} (or Corollary~\ref{iso_gr_sq_1}) for the exact description of the map.) \\
	
	But now consider the two-dimensional nonsingular quadratic form \(q \perp q \) over \(F\) given for all \((x_1,x_2) \in F \times F\) by 
	\[
	     (q \perp q)(x_1,x_2) = x_1^2 + x_2^2,
	\]
	which can be identified with the nondegenerate symmetric bilinear form \(b \perp b\) given for all \((x_1,x_2), (y_1,y_2) \in F \times F\) by
	\[
	   ( b \perp b)((x_1,x_2),(y_1,y_2)) = x_1y_1 + x_2y_2.
	\]
	Under the group isomorphism \(\partial^{1} \colon W_q(F)_{1} \to W(\overline{F}) \oplus W(\overline{F})\), the Witt class \({(q \perp q)}_W\) is sent to \(({(b \perp b)'}_W,0)\), where \((b \perp b)'\) is the nondegenerate two-dimensional symmetric bilinear form 
over \(\mathbb{F}_2\) given for all  \((x_1,x_2), (y_1,y_2) \in \mathbb{F}_2 \times \mathbb{F}_2 \) by
	\[
	        ( b \perp b)' ((x_1,x_2),(y_1,y_2)) = x_1y_1 + x_2y_2.
	\]
	The vector \((1,1)\) is clearly an isotropic vector for \( (b \perp b)'\), so that  \( (b \perp b)' \) is metabolic.\footnote{Note that it is a general fact that, when \(\charact F = 2\), the Witt group \(W(F)\) is a \(2\)-torsion group.} Hence \({(q \perp q)}_W \in \ker 
\partial^{1} = W_q(F)_{1/2}\). But it turns out that \({(q \perp q)}_W \not\in W_q(F)_{tame}\) (see Lemma~\ref{calculations_discrVal} \((e)\) and Corollary~\ref{Arason_2}). \qed
\end{example}

In this paper, each of the preceding homomorphisms (see \eqref{springer} and \eqref{filtration_isos}) will be obtained from a composition of two group homomorphisms. The first one is a canonical epimorphism (depending on $\varepsilon \in E$, but not on a 
choice of uniformizing parameters)
\begin{equation}\label{gr_iso}
  W_q(F)_{\varepsilon} \to W_g(F,\varepsilon),
\end{equation}
whose kernel is $W_q(F)_{< \varepsilon}$ when $\varepsilon > 0$. The image \(W_g(F,\varepsilon)\) of the map is a Witt-like group of the graded field $\gr_v(F)$ associated with the filtration of $F$ defined by $v$. That first canonical epimorphism is followed by 
an isomorphism which, in contrast, does depend on a choice of uniformizing parameters: for each $\varepsilon \in \frac 1 2 \Gamma_F$ such that $0\leq \varepsilon \leq v(2)$, we have
\begin{equation}\label{gr_res}
W_g(F,\varepsilon)\iso
 \begin{cases}
   \bigoplus_{\Gamma_F / 2 \Gamma_F} W_q(\overline{F}) & \text{if } \varepsilon = 0 \\ 
    \bigoplus_{\frac {\left| \Gamma_F / 2 \Gamma_F \right|} 2} W_{ssq}(\overline{F}) & \text{if } \varepsilon \notin \Gamma_F \\
    \bigoplus_{\Gamma_F / 2 \Gamma_F} W_{sq}(\overline{F}) & \text{if } \varepsilon \in \Gamma_F  \text{ and } \varepsilon \neq v(2) \text{ and } \varepsilon \neq 0 \\
    \bigoplus_{\Gamma_F / 2 \Gamma_F} W(\overline{F}) & \text{if } \varepsilon = v(2).
  \end{cases}
\end{equation}

Informally, the groups \(W_g(F,\varepsilon)\) can be thought of as groups containing all the residue forms ``in a canonical way'', whereas their descriptions as groups over the residue field depend on a choice of uniformizing parameters. \\

3) Suppose $\charact \overline{F} = 2$. Then there are group isomorphisms 
\[
  \begin{aligned}
    &W_{sq}(\overline{F}) \iso \overline{F} \wedge_{\overline{F}^2} \overline{F} &\text{ and } &W_{ssq}(\overline{F}) \iso \overline{F} \otimes_{\overline{F}^2} \overline{F}
  \end{aligned}
\] 
(in particular $W_{ssq}(\overline{F})$ is nontrivial). Here $\overline{F}^2 \subset \overline{F}$ is the subfield of squares in $\overline{F}$. \\

The structure of the paper is the following. In Section~\ref{section-1}, we introduce the basics about the Witt-like groups of graded fields. Section~\ref{section-2} is devoted to the proof of the isomorphisms  $W_{sq}(F) \iso F \wedge_{F^2} F$ and $W_{ssq}(F) 
\iso F \otimes_{F^2} F$ for an arbitrary field $F$ of characteristic $2$, which are precisely Theorem~\ref{structure_sq} and Theorem~\ref{structure_ssq}. In Section~\ref{section-3}, we establish isomorphisms between the Witt-like groups (dealt with in 
Section~\ref{section-1}) of a graded field $\mathsf{F}$ and Witt-like groups of $\mathsf{F}_0$, the subfield of $\mathsf{F}$ consisting of \(0\) and all the nonzero homogeneous elements of degree \(0\). These last isomorphisms depend on a choice of 
uniformizing parameters. See Corollary~\ref{iso_gr_sq_1} of Proposition~\ref{decompsitionCan_1}, and see also Proposition~\ref{iso_gr_sq_2}. In particular, the proof of the isomorphisms~\eqref{gr_res} follows directly from the results of this section. The 
construction of Arason's filtration and the canonical epimorphism~\eqref{gr_iso} is treated in Section~\ref{section-4}, see Theorem~\ref{filtration}. Finally, in the last section, we recover Arason's results in the case that $F$ is a discretely valued field, with our 
methods (Corollary~\ref{Arason_2}).

\section{Shifted quadratic spaces}
\label{section-1}

Let $\Gamma$ be a divisible torsion-free abelian group. A $\Gamma$-graded commutative ring in which every nonzero homogeneous element is invertible is called a \emph{$\Gamma$-graded field}, and a $\Gamma$-graded module over a \(\Gamma\)-graded 
field is called a \emph{$\Gamma$-graded vector space}. Since $\Gamma$ is torsion-free, $\Gamma$-graded fields are domains and $\Gamma$-graded vector spaces are free modules. The rank of a graded vector space is called its \emph{dimension}. For 
more information, see for example  \cite{MR3328410}*{\S 2.1}. In this section, ${\mathsf{F}}$ is a $\Gamma$-graded field and $\mathsf{V}$ denotes an arbitrary finite-dimensional $\Gamma$-graded ${\mathsf{F}}$-vector space. We write $
\Gamma_{\mathsf{F}} = \set{\gamma \in \Gamma \tq {\mathsf{F}}_{\gamma} \neq \set {0}}$ for the \emph{grade set} of \(\mathsf{F}\), where \(\mathsf{F}_{\gamma}\) denotes the additive abelian subgroup of \(\mathsf{F}\) consisting of \(0\) and all the nonzero 
homogeneous elements of degree \(\gamma\). A similar notation is used for \(\Gamma\)-graded modules.

For $\varepsilon \in \Gamma$, an \emph {$\varepsilon$-shifted graded bilinear form} on \( \mathsf{V}\) is an ${\mathsf{F}}$-bilinear map $b \colon \mathsf{V} \times \mathsf{V} \to {\mathsf{F}}$ such that for all $\gamma, \delta \in \Gamma$, we have 
$b(\mathsf{V}_{\gamma}, \mathsf{V}_{\delta}) \subset {\mathsf{F}}_{\gamma + \delta + \varepsilon}$. Such a form is called \emph{nondegenerate} if the only $x \in \mathsf{V}$ such that $b(x,y)=0$ for all $y \in \mathsf{V}$ is $x=0$. When $\varepsilon = 0$, we 
call $b$ a \emph{graded bilinear form}.

A \emph{graded quadratic form} on $\mathsf{V}$ is a map $q \colon \mathsf{V} \to {\mathsf{F}}$ satisfying the following conditions involving $q$ and its \emph{polar form} $b_q \colon \mathsf{V} \times \mathsf{V} \to {\mathsf{F}}$ defined for all $v,w \in 
\mathsf{V}$ by
\[
b_q(v,w) = q(v+w)-q(v)-q(w) :
\]
\listeN{
	\item $q(\alpha x) = \alpha^2 q(x)$ for all $x \in  \mathsf{V}$, for all $\alpha \in \mathsf{F}$
    \item $b_q$ is a graded bilinear form on $\mathsf{V}$
	\item $q(\mathsf{V}_{\gamma}) \subset {\mathsf{F}}_{2\gamma}$ for all $\gamma \in \Gamma$.
}

The graded quadratic form $q$ is called \emph{nonsingular} if its polar form $b_q$ is nondegenerate. It is called \emph{totally singular} when $q(v+w)= q(v) + q(w)$ for $v,w \in \mathsf{V}$.

\defn{Let $\varepsilon \in \Gamma$. An \emph{$\varepsilon$-shifted quadratic space} on \(\mathsf{V}\) is a $3$-tuple $(\mathsf{V}, q, b)$ where $q$ is a graded quadratic form and $b$ is a nondegenerate $\varepsilon$-shifted graded symmetric bilinear form 
such that there exists $\alpha \in {\mathsf{F}}$ satisfying $b_q = \alpha b$. We call $\alpha$ a \emph{multiplier} of $(\mathsf{V}, q, b)$. Note that $\alpha = 0$ is allowed.
}

The \emph{trivial} \(\varepsilon\)-shifted quadratic space over \(\mathsf{F}\), denoted by \(0\), is the $\varepsilon$-shifted quadratic space $(\mathsf{V}, q, b)$ such that \(\mathsf{V} = 0\).
Observe that if an \(\varepsilon\)-shifted quadratic space is not trivial, then the multiplier \(\alpha\) is uniquely determined since \(b\) is nondegenerate\footnote{Every $\alpha \in \mathsf{F}$ is a multiplier of the trivial $\varepsilon$-shifted quadratic space $0$.}. 
Observe also that if $\alpha \neq 0$, then $\alpha$ is homogeneous and $\varepsilon = - \deg(\alpha)$. In the rest of the text, we will refer to the existence of the multiplier as the ``compatibility condition''. Here are the three main classes of $\varepsilon$-shifted 
quadratic spaces we are interested in.

\begin{example}\label{classes}
(i) Given ${\mathsf{F}}$ and $\varepsilon \in \Gamma$, an $\varepsilon$-shifted quadratic space $(\mathsf{V},q,b)$ over ${\mathsf{F}}$ is called an \emph{\(\varepsilon\)-shifted (quadratic) space of type $I$} over \(\mathsf{F}\) if it satisfies the additional 
conditions $\varepsilon = 0$ and $b = b_q$ (i.e. $\alpha = 1$ is a multiplier), or if it is trivial.
(ii) Given ${\mathsf{F}}$ and $\varepsilon \in \Gamma$, an $\varepsilon$-shifted quadratic space $(\mathsf{V},q,b)$ over ${\mathsf{F}}$ is called an \emph{\(\varepsilon\)-shifted (quadratic) space of type $II$} over \(\mathsf{F}\) if $q$ is totally singular and $b$ 
is alternating. In this case, $\alpha = 0$ can always be chosen as a multiplier. Those spaces are also called \emph{$\varepsilon$-shifted graded symplectic quadratic spaces} over \(\mathsf{F}\).
(iii) Given ${\mathsf{F}}$ and $\varepsilon \in \Gamma_{\mathsf{F}}$ and a nonzero homogeneous $\tau \in {\mathsf{F}}_{\varepsilon}$, an $\varepsilon$-shifted quadratic space $(\mathsf{V},q,b)$ over ${\mathsf{F}}$ is called \emph{\(\varepsilon\)-shifted 
(quadratic) space of type $\tau-III$} over \(\mathsf{F}\), if we have $q(v) = \tau^{-1} b(v,v)$ for all $v \in \mathsf{V}$. Here $\alpha = 2 \tau^{-1}$ is always a multiplier since $b_q(v,w) = q(v+w) - q(v) -q(w) = 2 \tau^{-1}b(v,w)$  for all $v, w \in V$.
\end{example}

Note that an \(\varepsilon\)-shifted quadratic space can belong to more than one of the three classes of Example~\ref{classes}, and therefore can be of more than one \emph{type} \(T \in \set{I, II, \tau-III \, | \, \tau \in \mathsf{F}^{\times}_{\varepsilon}}\).

Two $\varepsilon$-shifted quadratic spaces $(\mathsf{V}, q, b)$ and $(\mathsf{V}', q', b')$ over \(\mathsf{F}\) with a common multiplier are \emph{isometric} if there is a (degree-preserving) graded \(\mathsf{F}\)-linear isomorphism $L \colon \mathsf{V} \to 
\mathsf{V}'$ which commutes with both $q$, $q'$ and $b$, $b'$. The \emph{orthogonal sum} of two $\varepsilon$-shifted quadratic spaces $(\mathsf{V}, q, b)$ and $(\mathsf{V}', q', b')$ with a common multiplier $\alpha$ is defined in the obvious way. It is an $
\varepsilon$-shifted quadratic space of multiplier $\alpha$ again. It is important to note that if $(\mathsf{V},q,b)$ is an $\varepsilon$-shifted quadratic space and if $\mathsf{V}_1$ and $\mathsf{V}_2$ are graded subspaces of $\mathsf{V}$ such that $\mathsf{V} 
= \mathsf{V}_1 \oplus \mathsf{V}_2$ and $b = b_{|\mathsf{V}_1} \perp b_{|\mathsf{V}_2}$, then we automatically get $q = q_{|\mathsf{V}_1} \perp q_{|\mathsf{V}_2}$ by the compatibility condition. Note also that the three classes of Example~\ref{classes} are 
closed under orthogonal sums.

An $\varepsilon$-shifted quadratic space $(\mathsf{V}, q, b)$ is \emph{isotropic} if there exists a nonzero $v\in \mathsf{V}$ such that $q(v)=0=b(v,v)$, and such a nonzero \(v\) is called an \emph{isotropic} vector (if there is no isotropic vector, the space $
(\mathsf{V}, q, b)$ is \emph{anisotropic}). Observe that an isotropic space has always a homogeneous isotropic vector. Indeed, since $\Gamma$ is torsion-free, we can give it a total order (see for example \cite{MR3328410}*{Remark 2.2}), turning it into a 
totally ordered abelian group. Therefore, for an isotropic vector $\sum_{\gamma \in \Gamma} v_{\gamma} $(with \(v_{\gamma} \in \mathsf{V}_{\gamma}\)), the $v_{\gamma'}$ component, where $\gamma'$ is the smallest $\gamma \in \Gamma$ such that 
$v_{\gamma} \neq 0$, is clearly isotropic. An $\varepsilon$-shifted quadratic space $(\mathsf{V}, q, b)$ is \emph{metabolic} if it contains a graded subspace $\mathsf{L} \subset \mathsf{V}$ of dimension $\dim \mathsf{L} = \frac 1 2\dim \mathsf{V}$ such that 
$q(\mathsf{L}) = \set{0} = b(\mathsf{L},\mathsf{L})$. Such a space is called a \emph{Lagrangian} of $(\mathsf{V}, q, b)$. The orthogonal sum of two metabolic spaces is also metabolic (a Lagrangian can be chosen as the direct sum of the Lagrangians). 
Besides, for every $\varepsilon$-shifted quadratic space $(\mathsf{V},q,b)$, the orthogonal sum $(\mathsf{V}, q, b) \perp -(\mathsf{V},q,b)$ is metabolic (since $\set{(v,v) \tq v \in \mathsf{V}}$ is a Lagrangian), where $-(\mathsf{V},q,b) := (\mathsf{V},-q,-b)$.

\begin{example}
(i) An $\varepsilon$-shifted quadratic space $(\mathsf{V},q, b)$ of type $I$ can be identified with a graded quadratic space $(\mathsf{V},q)$. This identification preserves orthogonal sum and isotropy. Moreover, metabolic graded spaces of type $I$ are 
metabolic graded quadratic spaces, which are in fact hyperbolic graded quadratic spaces.
(ii) An $\varepsilon$-shifted quadratic space $(\mathsf{V},q, b)$ of type $II$ is isotropic if and only if $q$ is isotropic.
(iii) An $\varepsilon$-shifted quadratic space $(\mathsf{V},q, b)$ of type $\tau-III$ can be identified with a nondegenerate graded symmetric bilinear space $(\mathsf{V}, \tau^{-1}b)$. This identification preserves orthogonal sum, isotropy and metabolicity.
\end{example}

\subsection{Normalisation}

\lem[Normalisation]{\label{normalisation_gr}Let $\varepsilon \in \Gamma$ and $b \colon \mathsf{V} \times \mathsf{V} \to {\mathsf{F}}$ a nondegenerate $\varepsilon$-shifted graded symmetric bilinear form. The space $(\mathsf{V},b)$ can be decomposed as 
\[
(\mathsf{V},b) = (\mathsf{V}_1, b_1) \perp \dots \perp (\mathsf{V}_n,b_n) \perp (\mathsf{W}, b_{|\mathsf{W}}),
\]
where $(\mathsf{V}_i,b_i)$ is a one-dimensional nondegenerate $\varepsilon$-shifted graded symmetric bilinear graded subspace of $(\mathsf{V},b)$ for all $i=1, \dots, n$ (with $n$ possibly equal to $0$) and $\mathsf{W} \subset \mathsf{V}$ is a graded 
subspace such that $ b_{|\mathsf{W}}$ is a nondegenerate $\varepsilon$-shifted graded alternating symmetric bilinear form. Consequently, $n = 0$ if $b$ is alternating, and $\mathsf{W} = \set{0}$ if $b$ is anisotropic. Moreover, $\mathsf{W}$ has a 
homogeneous basis $e_1, f_1, \dots, e_r, f_r$ for some $r \geq 0$ such that
 \listeN{
	\item $b(e_i,f_j) = \delta_{i,j}$ for all $i,j=1 \dots, r$ (in particular:  $\deg e_i + \deg f_i + \varepsilon = 0$ for $i=1 \dots, r$)
	\item $b(e_i,e_j) = b(f_i,f_j)= 0$ for all $i,j= 1, \dots, r$.
}}

We call such a basis for $(\mathsf{W},b_{|\mathsf{W}})$ a \emph{homogeneous symplectic basis}. It induces a decomposition into two-dimensional spaces:
\[
	(\mathsf{W}, b_{|\mathsf{W}}) = (\mathsf{W}'_1, b'_1) \perp \dots \perp (\mathsf{W}'_r, b'_r)
\]
where $\mathsf{W}'_i = \sev \langle e_i, f_i\rangle \subset \mathsf{V}$ and $b'_i = b_{|\mathsf{W}'_i}$ for all \(i = 1, \dots, n\).

\pr{The existence of the first decomposition is well-known for an (ungraded) symmetric bilinear space $(V',b')$ when $b'$ is nondegenerate. See for example \cite{MR2427530}*{Corollary 1.9}. The graded case is similar. The existence of a symplectic basis is 
well-known for a (ungraded) bilinear space $(W',b')$ when $b'$ is nondegenerate and alternating. See for example \cite{MR2427530}*{Proposition 1.8}. The graded case is similar.
}

Note that we can deduce from Lemma~\ref{normalisation_gr} a corresponding decomposition for an $\varepsilon$-shifted quadratic space $(\mathsf{V},q,b)$ by the compatibility condition.

\subsection{Witt decomposition}

\prop{\label{decomposition_gr}Let \(\varepsilon \in \Gamma\). Every $\varepsilon$-shifted quadratic space $(\mathsf{V}, q, b)$ can be decomposed as follows: $(\mathsf{V}, q, b) \simeq \varphi_{an} \perp \mu_1 \perp \dots \perp \mu_n$, where $\varphi_{an}$ 
is a uniquely determined (up to isometry) anisotropic $\varepsilon$-shifted quadratic space and $\mu_i$ is a metabolic two-dimensional $\varepsilon$-shifted quadratic space for all $i \in \set{1, \dots, n}$ (with $n$ possibly equal to $0$). Moreover, if $
(\mathsf{V}, q, b)$ is metabolic, then $\varphi_{an} = 0$.
}

\pr{This result is well-known for (ungraded) nondegenerate symmetric bilinear spaces $(V,b)$, see for example \cite{MR2427530}*{Theorem 1.27}. The case of graded nondegenerate symmetric bilinear spaces $(\mathsf{V}, b)$ is similar since an isotropic 
vector for $b$ can always be chosen homogeneous. Adapting the proof to $\varepsilon$-shifted quadratic spaces $(\mathsf{V}, q,b)$ is straightforward, because the orthogonality for $b$ implies the orthogonality for $b_q$ by the compatibility condition.
}

Note that the decomposition given in Proposition~\ref{decomposition_gr} preserves the type, i.e. if $\varphi$ is an $\varepsilon$-shifted quadratic space of type $T$ ($T \in \set{I,II, \tau-III \, | \, \tau \in \mathsf{F}^{\times}_{\varepsilon}}$), then $\varphi_{an}$ and 
all the $\mu_i$ (\(i = 1, \dots, n\)) are also of type $T$.

Mimicking the usual construction, for every \(\varepsilon \in \Gamma\) and every type $T \in \set{I,II, \tau-III \, | \, \tau \in \mathsf{F}^{\times}_{\varepsilon}}$, we may define a Witt equivalence of $\varepsilon$-shifted quadratic spaces of type \(T\) over a given 
graded field ${\mathsf{F}}$ and endow the set  $W_{T}^{\varepsilon}({\mathsf{F}})$ of Witt-equivalence classes of $\varepsilon$-shifted quadratic spaces of type $T$ over \(\mathsf{F}\) with a group structure using the orthogonal sum. Then each equivalence 
class is represented by a unique anisotropic space by Proposition~\ref{decomposition_gr}\footnote{Clearly, this construction works in general for $\varepsilon$-shifted quadratic spaces of the same multiplier $\alpha$.}.

After canonical identifications, we have $W_{I}^{0}({\mathsf{F}}) = W_{q}({\mathsf{F}})$, which is the quadratic Witt group of ${\mathsf{F}}$, and $W_{\tau-III}^{\varepsilon}({\mathsf{F}}) = W({\mathsf{F}})$, which is the Witt group of ${\mathsf{F}}$. Moreover, we 
also write $W_{sq}^{\varepsilon}({\mathsf{F}}) = W_{II}^{\varepsilon}({\mathsf{F}})$ for the Witt group of $\varepsilon$-shifted graded symplectic quadratic spaces.

%%%%%%%%%%%%%%%%%%%%%%%%%%%%%%%%%%%%%%%%%%%%%%%%%%%%%%%%%%%%%%%%%%%%%%%%%%%%%%%%
%%%%%%%%%%%%%%%%%%%%%%%%%%%%%%%%%%%%%%%%%%%%%%%%%%%%%%%%%%%%%%%%%%%%%%%%%%%%%%%%
%%%%%%%%%%%%%%%%%%%%%%%%%%%%%%%%%%%%%%%%%%%%%%%%%%%%%%%%%%%%%%%%%%%%%%%%%%%%%%%%

\section{Separated and nonseparated symplectic quadratic spaces}
\label{section-2}
In this section, $F$ is a field of characteristic $2$ and $V$ is an arbitrary finite-dimensional $F$-vector space.

\subsection{Symplectic quadratic spaces}

\defn{\label{def_ssq}A \emph{symplectic quadratic space over $F$} is a $0$-shifted graded symplectic quadratic space $(V,q,b)$ over $F$ considered as a $\Gamma$-graded field with $\Gamma = 0$.
}

In other words, $(V,q,b)$ is a symplectic quadratic space if $q$ is a totally singular quadratic form on $V$ and $b$ is a nondegenerate alternating bilinear form on V. As a special case of Lemma~\ref{normalisation_gr}, such a space always admits a symplectic 
basis $e_1,f_1, \dots, e_r, f_r$ for some $r \geq 0$. By Proposition~\ref{decomposition_gr}, we can form a Witt group of symplectic quadratic spaces over \(F\), which will be denoted by $W_{sq}(F)$.

\subsubsection{Structure}

For $\alpha, \alpha' \in F$, we denote by $\langle \alpha \, , \alpha' \rangle$ the two-dimensional symplectic quadratic space $(F\times F, q, b)$ where $q((1,0)) = \alpha$, $q((0,1)) = \alpha'$ and $b((1,0),(0,1)) = 1$. The basis \((1,0)\), \((0,1)\) is called the 
\emph{standard basis} of \(\langle \alpha \, , \alpha' \rangle\). Note that we use the notation \(\langle \phantom{x} \, , \, \phantom{x} \rangle\) instead of the usual square brackets notation \( [ \phantom{x} \, , \, \phantom{x} ]\) since this last notation will be used 
hereunder (see Corollary~\ref{W_q__asGenerated}) with another meaning.

\lem{\label{lem_sq}For all $\alpha, \alpha', \beta, \beta' \in F$, and for all $\xi \in F^{\times}$,
\listeN
{
  \item $\langle \alpha \, , \alpha \rangle$ is metabolic
  \item $\langle \alpha \, , \xi^2 \alpha' \rangle \simeq \langle \xi^2 \alpha \, , \alpha'  \rangle$
  \item $\langle \alpha \, , \alpha' \rangle \perp \langle \beta \, , \beta' \rangle \simeq \langle \alpha + \beta \, , \alpha' \rangle \perp \langle \beta \, , \alpha' + \beta' \rangle$.
}}

\pr{We show $(3)$. If $e_1, f_1$ is the standard basis of $\langle \alpha \, , \alpha' \rangle$ and $e_2, f_2$ the one of $\langle \beta \, , \beta' \rangle$, then $e_1 + e_2$ and $f_1$ spans a subspace $\langle \alpha + \beta \, , \alpha' \rangle$ with orthogonal 
complement the subspace $\langle \beta \, , \alpha' + \beta' \rangle$ spanned by $e_2$ and $f_1 + f_2$. The rest of the proof is left to the reader.
}

By Lemma~\ref{lem_sq}, the map $F \times F \to W_{sq}(F)$ which sends $(\alpha, \alpha')$ to $\langle \alpha \, , \alpha' \rangle$ induces a well-defined group homomorphism
\[
  \Phi \colon F \wedge_{F^2} F \to W_{sq}(F).
\]

\thm{\label{structure_sq}$\Phi$ is an isomorphism.
}

\pr{Since every symplectic quadratic space is a sum of two-dimensional spaces, $\Phi$ is surjective. So, we only have to prove that $\Phi$ is injective. We begin by proving that if $\alpha, \alpha', \beta, \beta' \in F$ are such that $\langle \alpha \, , \alpha' 
\rangle \simeq \langle \beta \, , \beta' \rangle$, then $\alpha \wedge \alpha' = \beta \wedge \beta'$. Let $e,f$ be the standard basis of $\langle \alpha \, , \alpha'  \rangle = (V, q, b)$. Given that $\langle \alpha \, , \alpha' \rangle \simeq \langle \beta \, , \beta' 
\rangle$, we can find vectors $x,y \in V$ such that
\[
  \begin{aligned}
    &q(x) = \beta, \qquad &b(x,y) = 1, \qquad &q(y) = \beta'.
  \end{aligned}
\]
Write $x = x_1 e + x_2 f$ and $y = y_1e + y_2 f$ with $x_1,x_2,y_1,y_2 \in F$. The previous relations give
\[
  \begin{aligned}
    &\alpha x_1^2 + \alpha' x_2^2 = \beta, \qquad &x_1y_2+x_2y_1 = 1,\qquad &\alpha y_1^2 + \alpha' y_2^2 = \beta'.
  \end{aligned}
\]
Therefore,
\[
  (\beta \wedge \beta') = (\alpha x_1^2 + \alpha' x_2^2) \wedge (\alpha y_1^2 + \alpha' y_2^2) = (x_1y_2+x_2y_1)^2 \alpha \wedge \alpha' = \alpha \wedge \alpha'. 
\]
Now, we show that if $\alpha_1, \alpha'_1, \dots, \alpha_n, \alpha'_n \in F$ are such that $\langle \alpha_1 \, , \alpha'_1 \rangle \perp \dots \perp \langle \alpha_n \, , \alpha'_n \rangle$ is isotropic, then we can rewrite $\alpha_1 \wedge \alpha'_1 + \dots + 
\alpha_n \wedge \alpha'_n$ as a sum of $n-1$ products. We prove this assertion by induction on $n$. Suppose that the statement is true for the sums of $n-1$ terms. Let $x \in V$ be an isotropic vector for $\langle \alpha_1 \, , \alpha'_1 \rangle \perp \dots \perp 
\langle \alpha_n \, , \alpha'_n \rangle = (V, q, b)$. We can assume that in the expression $x = x_1 + \dots + x_n$ with $x_i \in \langle \alpha_i \, , \alpha'_i \rangle$, each $x_i \neq 0$; otherwise we are done by the induction hypothesis. We may then find $y_i \in 
\langle \alpha_i \, , \alpha'_i \rangle$ such that $x_i, y_i$ is a symplectic basis. By putting $\beta_i = q(x_i)$ and $\beta'_i = q(y_i)$, we have $\langle \alpha_i \, , \alpha_i' \rangle \simeq \langle \beta_i \, , \beta_i' \rangle$, hence by the first part of the proof $
\alpha_i \wedge \alpha'_i = \beta_i \wedge \beta'_i$. Then,
\[
  \begin{aligned}
    \alpha_1 \wedge \alpha'_1 + \dots + \alpha_n \wedge \alpha'_n &= \beta_1 \wedge \beta'_1 + \dots + \beta_n \wedge \beta'_n \\
    &= \left( \beta_1 + \dots + \beta_n \right) \wedge \beta'_1 + \sum_{i = 2}^n \beta_i \wedge (\beta'_1 + \beta'_i).
  \end{aligned}
\]
Since $x$ is isotropic, $\beta_1 + \dots + \beta_n = 0$, which proves the assertion. To complete the proof, we show by induction on $n$ that if $\Phi\left(\sum_{i=1}^n \alpha_i \wedge \alpha'_i \right) = 0$, then $\sum_{i=1}^n \alpha_i \wedge \alpha'_i = 0$. This 
is clear in the case $n = 1$, for if $\langle \alpha \, , \alpha' \rangle$ is metabolic, then $\langle \alpha \, , \alpha' \rangle \simeq \langle \beta \, , 0 \rangle$ for some $\beta \in F$, thus by the first part of the proof $\alpha \wedge \alpha'  = \beta \wedge 0 = 0$. If 
$n > 1$ and $\Phi\left(\sum_{i=1}^n \alpha_i \wedge \alpha'_i \right) = 0$, then  $\langle \alpha_1 \, , \alpha_1 \rangle \perp \dots \perp \langle \alpha_n \, , \alpha'_n \rangle$ is metabolic, therefore isotropic and we can rewrite $\sum_{i=1}^n \alpha_i \wedge 
\alpha'_i$ as a sum of $n-1$ terms. By the induction  hypothesis, we have $\sum_{i=1}^n \alpha_i \wedge \alpha'_i = 0$.
}

The proof of Theorem~\ref{structure_sq} is inspired by \cite{RH-18-2006}*{page 4}.

\subsection{Separated symplectic quadratic spaces}

The dual of an arbitrary finite-dimensional \(F\)-vector space $E$ is denoted by $E^{*}$. If $L \colon E \to G$ is a linear map between two \(F\)-vector spaces $E$ and $G$, we denote by $L^* \colon G^* \to E^*$ the dual map of $L$. If $U \subset E$ is a 
subspace, we write $U^{o} = \set{\varphi \in E^* \tq \varphi(u) = 0 \text{ for all } u \in U}$ for the subspace in \(E^*\) which is orthogonal to \(U\). If $e_1, \dots, e_n$ is a basis of $E$, then the notation $e_1^*, \dots, e_n^*$ refers to its dual basis in $E^*$.

\defn{A \emph{separated symplectic quadratic space (over F)} is a $3$-tuple $(V, q, q')$ where $q$ is a totally singular quadratic form on $V$ and $q'$ is a totally singular quadratic form on $V^*$.
}

Two separated symplectic quadratic spaces $(V_1, q_1, q'_1)$ and $(V_2, q_2, q'_2)$ are \emph{isometric} if there exists an $F$-linear isomorphism $L \colon V_1 \to V_2$ such that $L$ is an isometry between $q_1$ and $q_2$, and $L^*$ is an isometry 
between $q'_2$ and $q'_1$. The \emph{orthogonal sum} of two spaces is defined in the obvious way (note that $(V_1 \oplus V_2)^*$ identifies with ${V_1}^* \oplus {V_2}^*$). A separated symplectic quadratic space $(V, q, q')$ is \emph{isotropic} if $q$ or $q'$ 
is isotropic, and it is \emph{metabolic} if it contains a subspace $U \subset V$ such that $q(U) = \set{0}= q'(U^o)$. Such a space \(U\) is called a \emph{Lagrangian}\footnote{There is no condition on the dimension of \(U\).} of $(V, q, q')$. The orthogonal sum of 
two metabolic spaces is also metabolic (the direct sum of the Lagrangians is a Lagrangian). Besides, for every separated symplectic quadratic space $(V, q, q')$, the orthogonal sum $(V, q, q') \perp (V, q, q')$ is metabolic since $\set{(v,v) \tq v \in V}$ is a 
Lagrangian. 

\subsubsection{Symplectic quadratic spaces and separated symplectic quadratic spaces}\label{functorSSQ-SQ}
In this paragraph, we define a functor $\mathcal{U}$ from separated symplectic quadratic spaces to  symplectic quadratic spaces.

Given a separated symplectic quadratic space $\Phi = (V, q, q')$, we put $\mathcal{U}(\Phi) = (V \oplus V^*, q \perp q', b)$ where $b \colon (V \oplus V^*) \times (V \oplus V^*) \to F$ is the nondegenerate alternating bilinear form defined by $b((v, \varphi), (w, 
\psi)) = \psi(v) - \varphi(w)$. That is indeed a symplectic quadratic space. Besides, if $L \colon \Phi \to \Psi$ is an isometry between two separated symplectic quadratic spaces $\Phi$ and $\Psi$, then $\mathcal{U}(L) := L \oplus L^*$ is an isometry $\mathcal{U}
(L) \colon \mathcal{U}(\Phi) \to \mathcal{U}(\Psi)$. It is also easy to see that $\mathcal{U}(\Phi \perp \Psi) \simeq \mathcal{U}(\Phi) \perp  \mathcal{U}(\Psi)$, and if $\Phi$ is metabolic then $\mathcal{U}(\Phi)$ is also metabolic, since if $L$ is a Lagrangian for $
\Phi$ then $L\oplus L^o$ is a Lagrangian for $\mathcal{U}(\Phi)$. Finally, if $x \in V$ (or $\varphi \in V^*$) is an isotropic vector for $\Phi$, then $(x,0)$ (or $(0, \varphi)$) is an isotropic vector for $\mathcal{U}(\Phi)$.

\subsubsection{Normalisation}
For $\alpha, \alpha' \in F$, we denote by $\langle \alpha \, | \, \alpha' \rangle$ the one-dimensional separated symplectic quadratic space $(F, q, q')$ where $q(1) = \alpha$ and $q'(1^*) = \alpha'$. The basis \(1\) is called the \emph{standard basis} of $\langle 
\alpha \, | \, \alpha' \rangle$.

\rem{For $\alpha, \alpha', \beta, \beta' \in F$, consider the two separated symplectic quadratic spaces  $\langle \alpha \, | \, \alpha'  \rangle$ and $\langle \beta \, | \, \beta'  \rangle$. Let $\langle \alpha \, | \, \alpha'  \rangle = (V_1, q_1, q_1')$. We have $\langle 
\alpha \, | \, \alpha'  \rangle \simeq \langle \beta \, | \, \beta'  \rangle$ if and only if there exists a nonzero vector $x \in V_1$ such that $q_1(x) = \beta$ and $q'_1(x^*) = \beta'$.
}

\lem[Normalisation]{\label{normalisation_ssq}Every separated symplectic quadratic space $(V,q,q')$ splits up into an orthogonal sum of one-dimensional separated symplectic quadratic spaces. More precisely, there exist $\alpha_1, \alpha'_1, \dots, \alpha_n, 
\alpha'_n \in F$ such that
\[
  (V, q, q') \simeq \langle \alpha_1 \, | \, \alpha'_1 \rangle \perp \dots \perp \langle \alpha_n \, | \, \alpha'_n \rangle
\]
}

\pr{The proof can be obtained by a variant of the proof of the existence of a symplectic basis in Lemma~\ref{normalisation_gr}. If $\dim V > 0$, take a nonzero $e \in V$ and a $\varphi \in V^*$ such that $\varphi (e) = 1$ and write $U = \sev \langle e \rangle 
\subset V$. Then decompose $(V, q, q') \simeq (U,q_{|U}, q'_{|U^*}) \perp (\ker \varphi, q_{|\ker \varphi}, q'_{|\left( \ker \varphi \right)^*})$, where $(U,q_{|U}, q'_{|U^*}) \simeq \langle \alpha \, | \, \alpha' \rangle$ for $\alpha = q_{|U}(e)$ and $\alpha' = q'_{|U^*}
(\varphi_{|U})$.
}

\subsubsection{Witt decomposition}

\prop{\label{decomposition_ssq}Every separated symplectic quadratic space $(V, q, q')$ can be decomposed as follows: $(V,q,q') \simeq \varphi_{an} \perp \mu_1 \perp \dots \perp \mu_n$, where $\varphi_{an}$ is a uniquely determined (up to isometry) 
anisotropic separated symplectic quadratic space and $\mu_i$ is a metabolic line for all $i \in \set{1, \dots, n}$ (with $n$ possibly equal to $0$). Moreover, if $(V, q, q')$ is metabolic, then $\varphi_{an} = 0$.
}

\pr{The proof is a refinement of the proof of Proposition~\ref{decomposition_gr} (when $\Gamma = 0$) taking into account the separation of  $\mathcal{U}(V, q, q')$ into the two particular spaces $V$ and $V^*$. 
}

As for symplectic quadratic spaces, we may define a Witt equivalence of separated symplectic quadratic spaces over a given field $F$ and endow the set $W_{ssq}(F)$ of Witt-equivalence classes with a group structure using the orthogonal sum. Then each 
equivalence class is represented by a unique anisotropic space by Proposition~\ref{decomposition_ssq}.

\subsubsection{Structure}
\lem{For all $\alpha, \alpha', \beta, \beta' \in F$, and for all $\xi \in F^{\times}$,
\listeN
{
  \item $\langle \alpha \, | \, \xi^2 \alpha' \rangle \simeq \langle \xi^2 \alpha \, | \, \alpha'  \rangle$
  \item $\langle \alpha \, | \, \alpha' \rangle \perp \langle \beta \, | \, \beta' \rangle \simeq \langle \alpha + \beta \, | \, \alpha' \rangle \perp \langle \beta \, | \, \alpha' + \beta' \rangle$.
}}

\pr{For $(2)$, note that if $e_1$ is the standard basis of $\langle \alpha \, | \, \alpha' \rangle$ and $e_2$ the one of $\langle \beta \, | \, \beta' \rangle$, then the basis $(e_1 + e_2, e_2)$, whose dual basis is $(e_1^*,e_1^* + e_2^*)$, gives the result. $(1)$ is left 
to the reader.
}

By the previous lemma, the map $F \times F \to W_{ssq}(F)$ which sends $(\alpha, \alpha')$ to $\langle \alpha \, | \, \alpha' \rangle$ induces a well-defined group homomorphism
\[
  \Phi \colon F \otimes_{F^2} F \to W_{ssq}(F).
\]

\thm{\label{structure_ssq}$\Phi$ is an isomorphism.
}

\pr{The proof uses the same ideas as the proof of Theorem~\ref{structure_sq}. For example, we prove that if $\alpha, \alpha', \beta, \beta' \in F$ are such that $\langle \alpha \, | \, \alpha' \rangle \simeq \langle \beta \, | \, \beta' \rangle$, then $\alpha \otimes 
\alpha' = \beta \otimes \beta'$. Let $e$ be the standard basis of $\langle \alpha \, | \, \alpha'  \rangle = (V, q, q')$. Given that $\langle \alpha \, | \, \alpha' \rangle \simeq \langle \beta \, | \, \beta' \rangle$, we can find a vector $x \in V$ such that
\[
  \begin{aligned}
    &q(x) = \beta &q'(x^*) = \beta'.
  \end{aligned}
\]
Let's write $x = \xi e$ and $e^* = \xi x^*$ with $x_1 \in F^{\times}$. The previous relations give
\[
  \begin{aligned}
    &\alpha \xi^2 = \beta & \alpha' = \xi^2 \beta'.
  \end{aligned}
\]
Therefore,
\[
  (\beta \otimes \beta') = (\xi^2 \alpha) \otimes \beta' = \alpha \otimes (\xi^2 \beta') = \alpha \otimes \alpha'. 
\]
Now, we show by induction on $n$ that if $\alpha_1, \alpha'_1, \dots, \alpha_n, \alpha'_n \in F$ are such that $\langle \alpha_1 \, | \, \alpha'_1 \rangle \perp \dots \perp \langle \alpha_n \, | \, \alpha'_n \rangle$ is isotropic, then we can rewrite $\alpha_1 \otimes 
\alpha'_1 + \dots + \alpha_n \otimes \alpha'_n$ as a sum of $n-1$ products. Suppose that the statement is true for the sums of $n-1$ terms. Let $x$ be an isotropic vector for $\langle \alpha_1 \, | \, \alpha'_1 \rangle \perp \dots \perp \langle \alpha_n \, | \, 
\alpha'_n \rangle = (V, q, q')$. Suppose $x \in V$ (the case $x \in V^*$ is similar). We can assume that in the expression $x = x_1 + \dots + x_n$ with $x_i \in \langle \alpha_i \, | \, \alpha'_i \rangle$, each $x_i \neq 0$. Therefore $x_1, \dots, x_n$ is a basis $V$. 
Let $x^*_1, \dots, x^*_n \in V^*$ be its dual basis. By putting $\beta_i = q(x_i)$ and $\beta'_i = q'(x^*_i)$, we have $\langle \alpha \, | \, \alpha' \rangle \simeq \langle \beta \, | \, \beta' \rangle$, hence by the first part of the proof $\alpha_i \otimes \alpha'_i = 
\beta_i \otimes \beta'_i$. Then, we conclude the induction as in Theorem~\ref{structure_sq}. The rest of the proof is left to the reader.
}

The functor $\mathcal{U}$ of \S~\ref{functorSSQ-SQ} induces a group homomorphism $W_{ssq}(F) \to W_{sq}(F)$ which fits in the following commutative diagram, where the left vertical arrow is canonical.
\[
\begin{tikzcd}
{F \otimes_{F^2} F} \arrow{r}{\sim} \arrow{d} & {W_{ssq}(F)} \arrow{d} \\
{F \wedge_{F^2} F} \arrow{r}{\sim} & {W_{sq}(F)}
\end{tikzcd}
\]

%%%%%%%%%%%%%%%%%%%%%%%%%%%%%%%%%%%%%%%%%%%%%%%%%%%%%%%%%%%%%%%%%%%%%%%%%%%%%%%%
%%%%%%%%%%%%%%%%%%%%%%%%%%%%%%%%%%%%%%%%%%%%%%%%%%%%%%%%%%%%%%%%%%%%%%%%%%%%%%%%
%%%%%%%%%%%%%%%%%%%%%%%%%%%%%%%%%%%%%%%%%%%%%%%%%%%%%%%%%%%%%%%%%%%%%%%%%%%%%%%%

\section{Structure of the Witt group of shifted spaces of type \(T\)}
\label{section-3}

This section is built upon the work of Elomary and Tignol (\cite{MR2836073}*{Section 2}) in which the case of graded quadratic spaces is already handled.

Throughout the section, we use the same notation as in Section~\ref{section-1}. In particular, \(\mathsf{F}\) is a \(\Gamma\)-graded field and ${\mathsf{V}}$ denotes an arbitrary finite-dimensional \(\Gamma\)-graded ${\mathsf{F}}$-vector space. The grade set 
of ${\mathsf{V}}$ is $\Gamma_{\mathsf{V}} = \set{\gamma \in \Gamma \tq {\mathsf{V}}_{\gamma} \neq 0}$. It is a union of cosets of $\Gamma_{\mathsf{F}}$. For $\gamma \in \Gamma$, we put $[\gamma] = \gamma + \Gamma_{\mathsf{F}} \in \Gamma / 
\Gamma_{\mathsf{F}}$ and we define
\[
	{\mathsf{V}}_{[\gamma]} = \bigoplus_{\delta \in \Gamma_{\mathsf{F}}} {\mathsf{V}}_{\gamma + \delta}.
\]
Clearly, ${\mathsf{V}}_{[\gamma]}$ does not depend on the representative $\gamma$ and is a graded subspace of ${\mathsf{V}}$ such that $\dim_{\mathsf{F}} {\mathsf{V}}_{[\gamma]} = \dim_{{\mathsf{F}}_0} {\mathsf{V}}_{\gamma + \delta}$ for all $\delta \in 
\Gamma_{\mathsf{F}}$. Moreover, we have
\[
	{\mathsf{V}} = \bigoplus_{\Lambda \in \Gamma_{\mathsf{V}} / \Gamma_{\mathsf{F}}} {\mathsf{V}}_{\Lambda},
\]
which is called \emph{the canonical decomposition of ${\mathsf{V}}$}.

For the rest of this section, we fix some $\varepsilon \in \Gamma$. By an \emph{$\varepsilon$-shifted space of type $T$} over \(\mathsf{F}\), we mean an $\varepsilon$-shifted quadratic space over \(\mathsf{F}\) of type $I$, $II$ or $\tau-III$ for some $\tau \in 
\mathsf{F}^{\times}_{\varepsilon}$. Note that if $({\mathsf{V}}, q, b)$ is an $\varepsilon$-shifted space of type $T$ and $v \in \mathsf{V}$, then $v$ is isotropic if and only if $q(v) = 0$, hence $({\mathsf{V}}, q, b)$ anisotropic implies $\Gamma_{\mathsf{V}} 
\subset \frac 1 2 \Gamma_{\mathsf{F}}$.

Note also that if $({\mathsf{V}},q,b)$ is an $\varepsilon$-shifted space of type $I$ or $\tau-III$, then we always have $\varepsilon \in \Gamma_{\mathsf{F}} \leq \frac 1 2 \Gamma_{\mathsf{F}}$ by definition.

\lem{\label{origin_of_the_cases}Suppose $({\mathsf{V}},q,b)$ is a nonzero anisotropic $\varepsilon$-shifted symplectic quadratic space. Let $e, f \in {\mathsf{V}}$ be homogeneous elements such that $b(e,f) = 1$. Then 
\listeN{
	\item $\varepsilon \in \frac 1 2 \Gamma_{\mathsf{F}}$
	\item $\deg e$ and $\deg f$ are in the same equivalence class in $ \frac 1 2 \Gamma_{\mathsf{F}} / \Gamma_{\mathsf{F}}$ if and only if $\varepsilon \in \Gamma_{\mathsf{F}}$.
}}

\pr{The lemma follows from the equation $\deg e + \deg f + \varepsilon = 0$, since $\deg e, \deg f \in \frac 1 2 \Gamma_{\mathsf{F}}$ by anisotropy of $({\mathsf{V}},q,b)$.
}

For $\gamma \in \Gamma$, we put $\overline{\gamma} = -(\gamma + \varepsilon)$.

\lem{\label{involution}There exists a well-defined involution $\Gamma / \Gamma_{\mathsf{F}} \to \Gamma /\Gamma_{\mathsf{F}} \colon [\gamma] \mapsto \overline{[\gamma]} := [\overline{\gamma}]$. Assume moreover $\varepsilon \in \frac 1 2 
\Gamma_{\mathsf{F}}$. Then the above involution restricts to a well-defined involution $\frac 1 2 \Gamma_{\mathsf{F}} / \Gamma_{\mathsf{F}} \to \frac 1 2 \Gamma_{\mathsf{F}} / \Gamma_{\mathsf{F}}$. If $\varepsilon \in \Gamma_{\mathsf{F}}$, then $
[\gamma] = \overline{[\gamma]}$ for all $\gamma \in \frac 1 2 \Gamma_{\mathsf{F}}$. If $\varepsilon \in \frac 1 2 \Gamma_{\mathsf{F}} \setminus \Gamma_{\mathsf{F}}$, then there is no \(\gamma \in \frac 1 2 \Gamma_{\mathsf{F}}\) such that $[\gamma] = 
\overline{[\gamma]}$. For $\varepsilon \in \frac 1 2 \Gamma_{\mathsf{F}}$, the above involution also restricts to a well-defined involution $( \Gamma / \Gamma_{\mathsf{F}} ) \setminus (  \frac 1 2 \Gamma_{\mathsf{F}} / \Gamma_{\mathsf{F}})  \to ( \Gamma / 
\Gamma_{\mathsf{F}} ) \setminus (  \frac 1 2 \Gamma_{\mathsf{F}} / \Gamma_{\mathsf{F}} )$, and this involution has no fixed points when $\varepsilon \in \Gamma_{\mathsf{F}}$.
}

\pr{The proof is left to the reader.}

Now suppose \(\varepsilon \in \frac 1 2 \Gamma_{\mathsf{F}}\). The involution of Lemma~\ref{involution} induces a group action $\Z / 2 \Z \to \aut(\Gamma / \Gamma_{\mathsf{F}})$ that preserves $\frac 1 2 \Gamma_{\mathsf{F}} / \Gamma_{\mathsf{F}}$ and its 
complement in $\Gamma / \Gamma_{\mathsf{F}}$. Let $\mathcal P$ be the set of orbits  under the action, which is a partition of $\Gamma / \Gamma_{\mathsf{F}}$. Similarly, let $\mPp$ the induced partition on $\frac 1 2 \Gamma_{\mathsf{F}} / 
\Gamma_{\mathsf{F}}$, and $\mPm$ the induced partition on $(\Gamma / \Gamma_{\mathsf{F}}) \setminus (\frac 1 2 \Gamma_{\mathsf{F}} / \Gamma_{\mathsf{F}})$. So we have $\mP = \mPp \sqcup \mPm$. We call $\mPp$ the \emph{principal-parts set of $
\mathcal P$} and $\mPm$ the \emph{metabolic-parts set of $\mathcal P$}. Now, for $P \in \mathcal P$ and $\left( {\mathsf{V}}, q, b\right)$ an $\varepsilon$-shifted space of type $T$, let ${\mathsf{V}}_P = \bigoplus_{\Lambda \in P} {\mathsf{V}}_{\Lambda}$, 
$q_P = q_{|{\mathsf{V}}_P}$, $b_P = b_{|{\mathsf{V}}_P}$ and $\Phi_P = ({\mathsf{V}}_P, q_P, b_P)$. The next lemma shows that $\Phi_P = ({\mathsf{V}}_P, q_P, b_P)$ is also an $\varepsilon$-shifted space of (the same) type $T$.

\lem{\label{nondegenerate}Let $b \colon {\mathsf{V}} \times {\mathsf{V}} \to {\mathsf{F}}$ be a nondegenerate $\varepsilon$-shifted graded symmetric bilinear form. For all $\gamma \in \Gamma$, ${\mathsf{V}}_{[\gamma]}^{\perp} = \bigoplus_{[\delta] \neq 
\overline{[\gamma]}} {\mathsf{V}}_{[\delta]}$, $\dim {\mathsf{V}}_{[\gamma]} = \dim {\mathsf{V}}_{\overline{[\gamma]}}$, and $b_{P}$ is nondegenerate, where $P \in \mP$ is the orbit of $[\gamma]$.
}

\pr{The proof uses the same ideas as in \cite{MR2330736}*{Proposition 1.1}, and is based on the fact that $b({\mathsf{V}}_{[\gamma]}, {\mathsf{V}}_{[\delta]}) = 0$ whenever $\gamma + \delta + \varepsilon \notin \Gamma_{\mathsf{F}}$.
}

Our next goal is to describe the Witt group $W_T^{\varepsilon}(\mathsf{F})$ in terms of a Witt group of $\mathsf{F}_0$. As suggested by Lemma~\ref{origin_of_the_cases}, we will distinguish two cases: when $\varepsilon \in \Gamma_F$ and when $
\varepsilon \in \frac 1 2 \Gamma_F \setminus \Gamma_F$.

\subsection{Case 1 : if \texorpdfstring{$\varepsilon \in \Gamma_{\mathsf{F}}$}{}}

\prop{\label{struct_decomp2}Suppose $\varepsilon \in \Gamma_{\mathsf{F}}$ and $T \in \set{I, II, \tau-III \, | \, \tau \in \mathsf{F}^{\times}_{\varepsilon}}$. Let $\Phi= \left( {\mathsf{V}}, q, b \right)$ be an $\varepsilon$-shifted space of type $T$ over \(\mathsf{F}\). 
The canonical decomposition of ${\mathsf{V}}$ yields a decomposition of $\Phi$ :
\[
\begin{aligned}
	\Phi = \left(\uperp_{P \in \mPp} \Phi_{P} \right) \perp \Psi & \text{ where } \Psi = \uperp_{P \in \mPm} \Phi_{P},
\end{aligned}
\]
where $\Psi$ and all the $\Phi_{P}$ (for \(P \in \mPp\)) are $\varepsilon$-shifted spaces of type $T$ over \(\mathsf{F}\). Moreover $\Psi$ is metabolic.
}

\pr{The proof is similar to the one of \cite{MR2836073}*{Proposition 1}, which deals with the case of nonsingular graded quadratic spaces; it uses the fact that $q({\mathsf{V}}_{[\gamma]}) = {0}$ when $\gamma \notin \frac 1 2 \Gamma_{\mathsf{F}}$.
}

Recall that \(\varepsilon \in \Gamma_{\mathsf{F}}\) is fixed, and fix also a type \(T \in \set{I,II, \tau-III \, | \, \tau \in \mathsf{F}^{\times}_{\varepsilon}}\) over \(\mathsf{F}\). By a \emph{choice of uniformizing parameters (for \(\varepsilon\) and \(T\))}, we mean a pair 
$\mathcal C = (\rho, \pi)$, where $\rho \in {\mathsf{F}}^{\times}_{\varepsilon}$ is a nonzero homogeneous element and $\pi \colon \mPp \to {\mathsf{F}}$ is a map such that, for all $P \in \mPp$, we have $\pi_{P} := \pi(P) \in {\mathsf{F}}^{\times}_{2\gamma}$ for 
some $\gamma \in \frac 1 2 \Gamma_{\mathsf{F}}$ (depending on \(P\)) such that $P = \set{[\gamma]}$. Moreover, we always assume (i.e., it is part of the definition) $\rho = 1$ if $T=I$, and $\rho = \tau$ if $T=\tau-III$. Now, given a choice of uniformizing 
parameters $\mathcal C = (\rho, \pi)$, for every $\varepsilon$-shifted space $\Phi = ({\mathsf{V}}, q, b)$  of type $T$, let for each $P \in \mPp$,

\[
\begin{aligned}
	&q_{\mathcal C, P} \colon {\mathsf{V}}_{\frac 1 2 \delta} \to {\mathsf{F}}_0 \colon x \mapsto\pi_{P}^{-1}q(x) \\
	&b_{\mathcal C, P} \colon {\mathsf{V}}_{\frac 1 2 \delta} \times {\mathsf{V}}_{\frac 1 2 \delta} \to {\mathsf{F}}_0 \colon 
	(x,y) \mapsto (\pi_{P}\rho)^{-1} b(x, y),
\end{aligned}
\]
where $\delta := \deg \pi_{P} \in \Gamma_{\mathsf{F}}$.

Consider the field \(\mathsf{F}_0\) as a graded field with trivial graduation.\footnote{\emph{Trivial graduation} means here \(\Delta\)-grading with \(\Delta = 0\), i.e. \(\Delta\) is the trivial group.} Let \(T' \in \set{I, II, 1-III}\) be the following type over \(\mathsf{F}
_0\) : if \(T = I\) then \(T' = I\), if \(T = II\) then \(T' = II\), and if \(T = \tau-III\) for some \(\tau \in \mathsf{F}^{\times}_{\varepsilon}\) then \(T' = 1-III\) (where \(1 \in \mathsf{F}^{\times}_0\) is the multiplicative unit).
Since the restriction of $b$ to ${\mathsf{V}}_{\mathcal C, P} := {\mathsf{V}}_{\frac 1 2 \delta}$ is nondegenerate (if $e \in {\mathsf{V}}_{\frac 1 2 \delta}$ and $f \in {\mathsf{V}}_{-\frac 1 2 \delta - \varepsilon}$ are such that $b(e,f) = 1$, then $b_{\mathcal C, P}(e,
\pi_{\delta}\rho f) = 1$), we have that $\Phi_{\mathcal C, P} = ({\mathsf{V}}_{\mathcal C, P}, q_{\mathcal C, P},b_{\mathcal C, P})$ is a \(0\)-shifted quadratic space of type \(T'\) on the trivially graded ${\mathsf{F}}_0$-vector space ${\mathsf{V}}_{\mathcal C, P}
$. Indeed, if \(T = I\) then, since in this case $b_{\mathcal C, P} = b_{q_{\mathcal C, P}}$, we clearly have that \(\Phi_{\mathcal C, P}\) is of type \(T' = I\), i.e. \(\Phi_{\mathcal C, P}\) can be identified with the nonsingular \(\mathsf{F}_0\)-quadratic space \(({\mathsf{V}}_{\mathcal C, P}, q_{\mathcal C, P})\). If \(T = II\) then we clearly have that \(\Phi_{\mathcal C, P}\) is of type \(T'= II\), i.e. \(\Phi_{\mathcal C, P}\)  is a symplectic \(\mathsf{F}_0\)-quadratic space. And if \(T = \tau-III\)  then, since in this case \(q_{\mathcal C, P}(v) = b_{\mathcal C, P}(v,v)\) for all $v \in {\mathsf{V}}_{\mathcal C, P}$, we clearly have that \(\Phi_{\mathcal C, P}\) is of type \(T' = 1-III\), i.e. \(\Phi_{\mathcal C, P}\) can be identified with the nondegenerate symmetric \(\mathsf{F}_0\)-
bilinear space \(({\mathsf{V}}_{\mathcal C, P},b_{\mathcal C, P})\). 

Recall that $W_{T}^{\varepsilon}({\mathsf{F}})$ is the Witt group of $\varepsilon$-shifted spaces of type $T$ over \(\mathsf{F}\) and let $W_{T'}({\mathsf{F}}_0)$ be the Witt group of \(0\)-shifted spaces of type \(T'\) over \(\mathsf{F}_0\). In other words, after 
canonical identifications, if \(T = I\) then \(W_{T'}({\mathsf{F}}_0) = W_I({\mathsf{F}}_0) \simeq W_q({\mathsf{F}}_0)\) is the Witt group of nonsingular quadratic spaces over ${\mathsf{F}}_0$, if \(T = II\) then \(W_{T'}({\mathsf{F}}_0) = W_{II}({\mathsf{F}}_0) 
\simeq W_{sq}({\mathsf{F}}_0)\) is the Witt group of symplectic quadratic spaces over \(\mathsf{F}_0\), and if \(T = \tau-III\) then \(W_{T'}({\mathsf{F}}_0) = W_{1-III}({\mathsf{F}}_0) \simeq W({\mathsf{F}}_0)\) is the Witt group of nondegenerate symmetric 
bilinear spaces over \(\mathsf{F}_0\).

Note that, since $\varepsilon \in \Gamma_{\mathsf{F}}$, it follows that $|\mPp| = |\frac 1 2 \Gamma_{\mathsf{F}} / \Gamma_{\mathsf{F}}| = |\Gamma_{\mathsf{F}} / 2 \Gamma_{\mathsf{F}}|$.

\prop{\label{decompsitionCan_1}Assume that $\varepsilon \in \Gamma_{\mathsf{F}}$ and $T \in \set{I, II, \tau-III \, | \, \tau \in \mathsf{F}^{\times}_{\varepsilon}}$. Then for each choice of uniformizing parameters $\mathcal C$, the map that carries each $
\varepsilon$-shifted space $\Phi$ of type $T$ over \(\mathsf{F}\) to the collection $\left(\Phi_{\mathcal C, P}\right)_{P \in \mPp}$ induces a group homomorphism
\[
   W_{T}^{\varepsilon}({\mathsf{F}}) \iso \bigoplus_{\mPp} W_{T'}({\mathsf{F}}_0).
\]
That isomorphism depends on the choice of the uniformizing parameters.
}

\pr{The proof is routine. See \cite{MR2836073}*{Proposition 2} for the case of nonsingular graded quadratic forms. The general case is similar.
}

\cor{Assume that $\varepsilon \in \Gamma_{\mathsf{F}}$. Then for each choice of uniformizing parameters $\mathcal C$, 
\begin{enumerate}
\label{iso_gr_sq_1} 
\item the map that carries each $\varepsilon$-shifted graded symplectic quadratic space $\Phi$ to the collection $\left(\Phi_{\mathcal C, P}\right)_{P \in \mPp}$ induces a group homomorphism
\[
	W_{sq}^{\varepsilon}({\mathsf{F}}) \iso \bigoplus_{\mPp} W_{sq}({\mathsf{F}}_0).
\]
\item the map that carries each nonsingular graded quadratic space $\Phi$ to the collection $\left(\Phi_{\mathcal C, P}\right)_{P \in \mPp}$ induces a group homomorphism
\[
	W_{q}({\mathsf{F}}) \iso \bigoplus_{\mPp} W_{q}({\mathsf{F}}_0).
\]
\item the map that carries each nondegenerate graded symmetric bilinear space $\Phi$ to the collection $\left(\Phi_{\mathcal C, P}\right)_{P \in \mPp}$ induces a group homomorphism
\[
	W({\mathsf{F}}) \iso \bigoplus_{\mPp} W({\mathsf{F}}_0).
\]
\end{enumerate}
Those isomorphisms depend on the choice of the uniformizing parameters.
}

\subsection{Case 2 : if \texorpdfstring{$\varepsilon \in \frac 1 2 \Gamma_{\mathsf{F}} \setminus \Gamma_{\mathsf{F}}$}{}}

In this case, Lemma~\ref{involution} shows that each orbit in $\mPp$ has exactly two different elements of $\frac 1 2 \Gamma_{\mathsf{F}} / \Gamma_{\mathsf{F}}$. This case occurs only for $\varepsilon$-shifted space of type $II$.
\prop{\label{struct_decomp}Suppose $\varepsilon \in \frac 1 2 \Gamma_{\mathsf{F}} \setminus \Gamma_{\mathsf{F}}$. Let $\Phi= \left( {\mathsf{V}}, q, b \right)$ be an $\varepsilon$-shifted graded symplectic quadratic space over \(\mathsf{F}\). 
The canonical decomposition yields a decomposition of $\Phi$ : 
\[
\begin{aligned}
	\Phi = \left(\uperp_{P \in \mPp} \Phi_{P} \right) \perp \Psi & \text{ where } \Psi = \uperp_{P \in \mPm} \Phi_{P},
\end{aligned}
\]
For each $P = \set{\Lambda, \overline{\Lambda}} \in \mPp$, $\Phi_{P} = ({\mathsf{V}}_P, q_P, b_P)$ is an $\varepsilon$-shifted graded symplectic quadratic space over \(\mathsf{F}\), and both ${\mathsf{V}}_{\Lambda}$ and ${\mathsf{V}}_{\overline{\Lambda}}$ 
are totally isotropic subspaces for $b_{P}$. Moreover, $\Psi$ is a metabolic $\varepsilon$-shifted graded symplectic quadratic space over \(\mathsf{F}\).
}

\pr{The existence of the decomposition is clear by Lemma~\ref{nondegenerate}. Now, let $\gamma \in \Gamma$. If $\gamma \in \frac 1 2 \Gamma_{\mathsf{F}}$, then $\gamma + \gamma + \varepsilon \notin \Gamma_{\mathsf{F}}$ since $\varepsilon \notin 
\Gamma_{{\mathsf{F}}}$, hence ${\mathsf{V}}_{[\gamma]}$ is a totally isotropic subspace for $b$. To complete the proof, it remains only to prove the metabolicity of $\Phi_{P}$ when $P \in \mPm$. But this is clear, since $q({\mathsf{V}}_{P}) = 0$ when $P \in 
\mPm$, and since $b_P$ is alternating.
}

Recall that \(\varepsilon \in \frac 1 2 \Gamma_{\mathsf{F}} \setminus \Gamma_{\mathsf{F}}\) is fixed. For each $\varepsilon$-shifted graded symplectic quadratic space $\Phi = ({\mathsf{V}}, q, b)$ and each $\gamma \in \frac 1 2 \Gamma_{\mathsf{F}}$, 
consider the ${{\mathsf{F}}}_0$-bilinear map $b_{|{\mathsf{V}}_{\overline{\gamma}} \times {\mathsf{V}}_{\gamma}} \colon {\mathsf{V}}_{\overline{\gamma}} \times {\mathsf{V}}_{\gamma} \to {\mathsf{F}}_{0} \colon (x,y) \mapsto b(x,y)$. This map is 
nondegenerate because for every nonzero $x \in {\mathsf{V}}_{\overline{\gamma}}$, there exists a homogeneous $y \in {\mathsf{V}}$ such that $b(x,y)=1$, that is $\deg y = -\overline{\gamma} - \varepsilon = \gamma$. Since moreover $\dim_{{\mathsf{F}}_0}
{\mathsf{V}}_{\gamma} = \dim_{\mathsf{F}} {\mathsf{V}}_{[\gamma]} = \dim_{\mathsf{F}} {\mathsf{V}}_{\overline{[\gamma]}} = \dim_{{{\mathsf{F}}}_0} {\mathsf{V}}_{\overline{\gamma}}$, where the middle equation comes from Lemma~\ref{nondegenerate}, we 
obtain a linear isomorphism $\widehat{b}_{\gamma} \colon {\mathsf{V}}_{\overline{\gamma}} \to {\mathsf{V}}_{\gamma}^* \colon z \mapsto b_{|{\mathsf{V}}_{\overline{\gamma}} \times {\mathsf{V}}_{\gamma}}(z, \cdot)$.

By a \emph{choice of uniformizing parameters (for \(\varepsilon\))}, we mean a pair $\mathcal C = (\rho, \pi)$, where $\rho \in {\mathsf{F}}_{2 \varepsilon}$ is a nonzero homogeneous element and $\pi \colon{\mathcal P}_{p} \to {\mathsf{F}}$ is a map such that 
for all $P \in {\mathcal P}_{p}$, we have $\pi_P := \pi(P) \in {{\mathsf{F}}}^{\times}_{2\gamma}$ for some $\gamma \in \frac 1 2 \Gamma_{\mathsf{F}}$ (depending on \(P\)) such that $P = \set{[\gamma], [\overline{\gamma}]}$. Now, given a choice of uniformizing 
parameters $\mathcal C = (\rho, \pi)$, for every $\varepsilon$-shifted graded symplectic quadratic space $\Phi = ({\mathsf{V}}, q, b)$, let for each $P \in {\mathcal P}_{p}$,
\[
\begin{aligned}
	&q_{\mathcal C, P} \colon {\mathsf{V}}_{\frac 1 2 \delta} \to {\mathsf{F}}_0 \colon x \mapsto \pi_P^{-1}q(x) \\
	&q'_{\mathcal C, P} \colon {\mathsf{V}}_{\frac 1 2 \delta}^* \to {\mathsf{F}}_0 \colon \widehat{b}_{\frac 1 2 \delta}(z) \mapsto \pi_P\rho q(z),
\end{aligned}
\]
where $\delta := \deg \pi_P \in \Gamma_{\mathsf{F}}$.

It is clear that $\Phi_{\mathcal C, P} = ({\mathsf{V}}_{\frac 1 2 \delta}, q_{\mathcal C, P}, q'_{\mathcal C, P})$ is a separated symplectic quadratic space over ${\mathsf{V}}_{\frac 1 2 \delta}$ viewed as an ${\mathsf{F}}_0$-vector space.

\prop{\label{iso_gr_sq_2}Assume that $\varepsilon \in \frac 1 2 \Gamma_{\mathsf{F}} \setminus \Gamma_{\mathsf{F}}$. Then for each choice of uniformizing parameters $\mathcal C$, the map that carries each $\varepsilon$-shifted graded symplectic 
quadratic space $\Phi$ to the collection $\left(\Phi_{\mathcal C, P}\right)_{P \in \mPp}$ induces a group homomorphism
\[
	W_{sq}^{\varepsilon}({\mathsf{F}}) \iso  \bigoplus_{\mPp} W_{ssq}({\mathsf{F}}_0).
\]
That isomorphism depends on the choice of the uniformizing parameters.
}

\pr{The proof is routine and is based on the same ideas as the ones of the proof of Proposition~\ref{struct_decomp2}, except that we use the nondegenerate alternating bilinear form $b$ of $\Phi$ to identify $V_{\frac 1 2  \overline{\delta_P}}$ with $V^*_{\frac 1 
2 \delta_P}$ for each $P \in \mPp$, where $\delta_P = \deg (\pi_P)$ if we let $\mathcal C = (\rho, \pi)$.
}

%%%%%%%%%%%%%%%%%%%%%%%%%%%%%%%%%%%%%%%%%%%%%%%%%%%%%%%%%%%%%%%%%%%%%%%%%%%%%%%%
%%%%%%%%%%%%%%%%%%%%%%%%%%%%%%%%%%%%%%%%%%%%%%%%%%%%%%%%%%%%%%%%%%%%%%%%%%%%%%%% 
%%%%%%%%%%%%%%%%%%%%%%%%%%%%%%%%%%%%%%%%%%%%%%%%%%%%%%%%%%%%%%%%%%%%%%%%%%%%%%%%

\section{Arason's filtration}
\label{section-4}

Let $F$ be a field of arbitrary characteristic and let $v \colon F \to \Gamma \cup \set{\infty}$ be a valuation, where $\Gamma$ is an arbitrary totally ordered (hence torsion-free, see \cite{MR2215492}*{Lemma 2.1.1}) abelian group. Without loss of generality, we 
may also assume $\Gamma$ divisible, since we may substitute for $\Gamma$ its divisible hull (cf. \cite{MR2215492}*{Proposition 1.2.4 \& Proposition 2.1.2}). Denote by $\overline{F}$ the residue field. Let also $V$ be a finite-dimensional $F$-vector space. 
We recall from \cite{MR3328410}*{\S 3.1.1} that a \emph{$v$-value function} is a map $\alpha \colon V \to \Gamma \cup \set{\infty}$ such that for all $x,y \in V$ and $\lambda \in F$
\begin{itemize}
\item[(i)] 
  $\alpha(x) = 0$ if and only if $x = 0$
\item[(ii)]
  $\alpha(\lambda x) = v(\lambda) + \alpha(x)$
\item[(iii)]
   $\alpha(x + y) \geq \min \set{\alpha(x), \alpha(y)}$.
\end{itemize}

The $v$-value function is called a \emph{$v$-norm} if there is a basis $(e_i)_{i=1,\dots, n}$ of $V$ such that 
\[
	\alpha \left(\sum_{i=1}^n \lambda_i e_i \right) = \min \set{\alpha(\lambda_i e_i) \tq i=1, \dots, n} \qquad \text{for all } \lambda_1, \dots, \lambda_n \in F.
\]
Such a basis is called a \emph{splitting basis} for $\alpha$. For example, it turns out that if $F$ is maximally complete for $v$ (e.g., $F$ is complete and $v$ is discrete), then every $v$-value function on $V$ is a $v$-norm (see \cite{MR3328410}*{Proposition 
3.8}). Nevertheless, maximal completeness is not required for the results of this section.

The value function $\alpha$ yields a filtration of $V$. We also recall from \cite{MR3328410}*{\S 3.1.1} the construction of \(\gr_{\alpha}(V)\), the \emph{associated graded vector space} over $\gr_v(F)$. For each \(\gamma \in \Gamma\) we let
\begin{align*}
   & V_{\geq \gamma} = \set{x \in V \tq \alpha(x) \geq \gamma}, & V_{> \gamma} = \set{x \in V \tq \alpha(x) > \gamma}, \quad &V_{\gamma} = V_{\geq \gamma} / V_{> \gamma},    
\end{align*}
and we define
\[
   \gr_{\alpha}(V) = \bigoplus_{\gamma \in \Gamma} V_{\gamma}.
\]
Similarly, let \(\gr_v(F)\) be the \(\Gamma\)-graded ring associated with the filtration of \(F\) defined by the valuation \(v\). The field structure on $F$ induces canonically a structure of graded ring on $\gr_v(F)$, for which every nonzero homogeneous element is 
invertible. Similarly, the $F$-vector space structure on $V$ induces a structure of $\gr_v(F)$-module on $\gr_{\alpha}(V)$. In particular, every $V_{\gamma}$ is a $F_0$-vector space. For each nonzero $x \in V^{\times}$, we let
\[
	\widetilde{x} = x + V_{>\alpha(x)} \in V_{\alpha(x)} \subset \gr_{\alpha}(F).
\]
We also set $\widetilde{0} = 0$ and use a similar notation for elements in \(\gr_v(F)\). Note that $F_0 = \overline{F}$ , so that $\widetilde{x} = \overline{x}$ if $v(x)=0$. It is shown in \cite{MR3328410}*{Corollary 3.6} that a family of vectors $(e_i)_{i=1,\dots,n}$ 
of $V$ is a splitting basis for $\alpha$ if and only if $(\widetilde{e_i})_{i=1,\dots,n}$ is a homogeneous $\gr_{v}(F)$-basis of $\gr_{\alpha}(V)$, and that $\alpha$ is a \(v\)-norm if and only if $\dim_{\gr_v(F)} \gr_{\alpha}(V) = \dim_F V$. We also write $
\Gamma_F : = \Gamma_{\gr_v(F)}$ for \emph{the value group of $v$}.

\subsection{Depth of norms and induced spaces}

If $q \colon V \to F$ is a quadratic form, the \emph{polar form of $q$} is the symmetric bilinear form $b_{q} \colon V \times V \to F$ defined for all $v, w \in V$ by 
\[
  b_{q}(v,w) = q(v + w) - q(v) - q(w).
\]

\defn{\label{compatibility}Let $q \colon V \to F$ be a quadratic form, with polar form $b \colon V \times V \to F$ and let $\varepsilon \in \Gamma$, $\varepsilon \geq 0$. We say that a $v$-norm $\alpha \colon V \to \Gamma \cup \set{\infty}$ is \emph{compatible 
of depth $\varepsilon$} (or \emph{\(\varepsilon\)-compatible}) \emph{with $q$} if
\begin{enumerate}
  \item[(a)] $v(b(x,y)) \geq \alpha(x) + \alpha(y) + \varepsilon$  for all $x, y \in V$
  \item[(b)] $v(q(x)) \geq 2 \alpha(x)$ for all $x \in V$
  \item[(c)] for all nonzero $x \in V$, there exists a nonzero $y \in V$ such that $v(b(x,y)) = \alpha(x) + \alpha(y) + \varepsilon$.
\end{enumerate}
Such a compatible $v$-norm is said to be \emph{tame} if $\varepsilon = 0$.
}

Note that it suffices to check conditions $\mathrm{(a)}$ and $\mathrm{(b)}$ for a splitting basis of $\alpha$.
Observe also that if the quadratic form \(q\) admits a \(\varepsilon\)-compatible \(v\)-norm (for some \(\varepsilon\)), then \(q\) is nonsingular. 
Finally note that Elomary and Tignol consider tame compatible \(v\)-norms in \cite{MR2836073}*{Section 3}.

If $\varphi = (V,q)$ is a quadratic space over $F$ with polar form $b$ and if $\alpha$ is a $v$-norm on $V$ which is \(\varepsilon\)-compatible with $q$ for some depth $\varepsilon \geq 0$, we set for all nonzero \(x,y \in V\),
\[
	\widetilde{b}_{\alpha}(\widetilde{x}, \widetilde{y}) = b(x,y) + F_{>\alpha(x) + \alpha(y) + \varepsilon} \in \gr_{v}(F)_{\alpha(x) + \alpha(y) + \varepsilon},
\]
and we extend \(\widetilde{b}_{\alpha}\) to a map 
\[
   \widetilde{b}_{\alpha} \colon \gr_{\alpha}(V) \times \gr_{\alpha}(V) \to \gr_{v}(F)
\]
by bilinearity. This map is an $\varepsilon$-shifted graded symmetric bilinear form.
Condition $c$ in Definition~\ref{compatibility} means exactly that $\widetilde{b}_{\alpha}$ is nondegenerate. We may also define a graded quadratic form 
\[
   \widetilde{q}_{\alpha} \colon \gr_{\alpha}(V) \to \gr_{v}(F)
\]
which satisfies for all nonzero $x \in V$,
\[
	\widetilde{q}_{\alpha}(\widetilde{x}) =  q(x) + F_{> 2\alpha(x)} \in \gr_{v}(F)_{2\alpha(x)}.
\]
Indeed, if $\varepsilon > 0$, the above formula can be extended to define a totally singular quadratic form on $\gr_{\alpha}(V)$; whereas if $\varepsilon = 0$, for $\widetilde{x}_{\gamma}\in V_{\gamma}$, we extend the above formula by setting 
\[
   \widetilde{q}_{\alpha}\left(\sum_{\gamma \in \Gamma} \widetilde{x}_{\gamma}\right) 
   = \sum_{\gamma \in \Gamma} \widetilde{q}_{\alpha}\left(\widetilde{x}_{\gamma}\right) 
      + \sum_{\gamma < \delta} \widetilde{b}_{\alpha}\left(\widetilde{x}_{\gamma}, \widetilde{x}_{\delta}\right),
\]
which defines a nonsingular quadratic form of polar form $\widetilde{b}_{\alpha}$. The straightforward verifications are omitted. We put 
\[
   \widetilde{\varphi}_{\alpha} = (\gr_{\alpha}(V), \widetilde{q}_{\alpha}, \widetilde{b}_{\alpha})
\] 
and we call $\widetilde{\varphi}_{\alpha}$ the \emph{$\varepsilon$-shifted (quadratic) space induced by $\varphi$ (and $\alpha$)}. Note that if $\varepsilon < v(2)$ then $\widetilde{b}_{\alpha}$ is alternating since for all nonzero $x \in V$, 
\[
    v(b(x,x)) = v(q(x)) + v(2) > 2 \alpha(x) + \varepsilon.
\] 
We also have that if $\varepsilon = v(2) < \infty$, then $\widetilde 2$ is a nonzero homogeneous element (but we have $\widetilde 2 = \widetilde 1 + \widetilde 1$ only if $v(2) = 0$) and $\widetilde{q}_{\alpha}(v) =  \widetilde{2}^{-1} \widetilde{b}_{\alpha}(v,v)$ for 
all $v \in \gr_{\alpha}(V)$. Therefore, for $\varepsilon \in \Gamma$ and $\varepsilon \geq 0$, if $\varepsilon = 0$ then $\widetilde{\varphi}_{\alpha}$, which is an $\varepsilon$-shifted space of type $I$ over \(\gr_v(F)\), can be identified with the usual induced 
graded quadratic space over $(\gr_{\alpha}(V), \widetilde{q}_{\alpha})$ over $\gr_v(F)$ (the ones considered in \cite{MR2836073}*{Section 3}). If $0 < \varepsilon < v(2)$, then $\widetilde{\varphi}_{\alpha}$ is an $\varepsilon$-shifted graded symplectic 
quadratic space \(\gr_v(F)\). Finally, if $\varepsilon = v(2)$ then $\widetilde{\varphi}_{\alpha}$, which is an $\varepsilon$-shifted space of type $\widetilde{2}-III$ over \(\gr_v(F)\), can be identified with the nondegenerate graded symmetric bilinear space $
(\gr_{\alpha}(V), \widetilde{2}^{-1}\widetilde{b}_{\alpha})$.

\subsection{The filtration}

In this section, we construct and describe Arason's filtration.

\lem{\label{sum_of_norms}Let $\varphi_1 = (V_1, q_1)$ and $\varphi_2 = (V_2, q_2)$ be two quadratic spaces over $F$, and let $\alpha_1, \alpha_2$ be $v$-norms on $V_1$, $V_2$ that are \(\varepsilon\)-compatible with $q_1$ and $q_2$ respectively, for 
some depth $\varepsilon \geq 0$. Define $\alpha_1 \oplus \alpha_2 \colon V_1 \oplus V_2 \to \Gamma \cup \set{\infty}$ by 
\[
\begin{aligned}
	&(\alpha_1 \oplus \alpha_2)(x_1, x_2) = \min \set{\alpha(x_1), \alpha(x_2)} &\text{ for $x_1 \in V_1$ and $x_2 \in V_2$.}
\end{aligned}
\]
Then $\alpha_1 \oplus \alpha_2$ is a $v$-norm on $V_1 \oplus V_2$ which is \(\varepsilon\)-compatible with $q_1 \oplus q_2$, and there is a canonical identification of \(\varepsilon\)-shifted quadratic spaces
\[
	\widetilde{\left( \varphi_1 \perp \varphi_2 \right)}_{\alpha_1 \oplus \alpha_2} = \widetilde{\varphi_1}_{\alpha_1} \perp \widetilde{\varphi_2}_{\alpha_2}
\]
}

\pr{The case of tame $v$-norms is treated in \cite{MR2836073}*{Lemma 6}, the general case is similar.
}

\lem{\label{depth_of_norms}If a quadratic space $\varphi = (V,q,b)$ admits a compatible $v$-norm $\beta'$ of depth $\delta \geq 0$, then it also admits a compatible \(v\)-norm $\beta$ of depth $\gamma$ for all $\gamma > \delta$. 
}

\pr{Suppose that $\beta'$ is a compatible $v$-norm of depth $\delta \geq 0$. Then we construct $\beta$ by letting $\beta = \beta' - \frac 1 2 (\gamma - \delta)$, so that for all $x,y \in V$
\[
\begin{aligned}
	&\beta'(x) + \beta'(y) + \delta = \beta(x) + \beta(y) + \gamma &\text{ and } &\beta(x) < \beta'(x).
\end{aligned}
\]
Then $\beta$ is a $v$-norm which is \(\gamma\)-compatible with $\varphi$.
}

\lem{\label{existence_of_norms_lemma} ~
\listeN{
  \item Let \(q \colon F \times F \to F \colon (x,y) \mapsto xy\) be a two-dimensional hyperbolic quadratic form. Then \(q\) admits a tame compatible \(v\)-norm.
  \item Suppose $v(2) > 0$ and let \[q \colon F \times F \to F \colon (x,y) \mapsto ax^2  +xy + by^2\] (with $a, b \in F$) be a two-dimensional quadratic form such that $v(4ab) > 0$. Then $q$ is nonsingular. \\
  If $v(a) + v(b) \leq 0$, then $q$ admits a compatible $v$-norm of depth \[\varepsilon = - \frac 1 2 (v(a) + v(b)) \in \frac 1 2 \Gamma_F\] satisfying  $0 \leq \varepsilon < v(2)$. \\
  If $v(a) + v(b) > 0$, then $q$ admits a compatible $v$-norm of depth $\varepsilon = 0$.
  \item Suppose $\charact(F) = 0$  and let \[q \colon F \to F \colon x \mapsto ax^2\] (with $a \in F^{\times}$) be a one-dimensional nonsingular quadratic space. Then $q$ admits a compatible $v$-norm of depth $\varepsilon = v(2)$.
}}

\pr{(1) It is easy to check that the map $\alpha \colon F \times F \to \Gamma \cup \set{\infty}$ defined for all $x,y \in F$ by \[\alpha(x,y) = \min\set{v(x), v(y)}\] is a tame compatible $v$-norm with \(q\). \\
(2) Note that the condition $v(4ab) > 0$ implies that $q$ is nonsingular. \\
If $v(a) + v(b) \leq 0$, then \[0 \leq - \frac 1 2 (v(a) + v(b)) < v(2)\] and the map $\alpha \colon F \times F \to \Gamma \cup \set{\infty}$ defined for all $x,y \in F$ by \[\alpha(x,y) = \min\set{v(x) + \frac 12 v(a), v(y) + \frac 1 2 v(b)}\] is a $v$-norm which is 
compatible of depth \[\varepsilon = - \frac 1 2 (v(a) + v(b))\] with $q$. \\
If $v(a) + v(b) > 0$ and for example $v(b) < 0$, then $v(a) > -v(b)$ and the map $\alpha \colon F \times F \to \Gamma \cup \set{\infty}$ defined for all $x,y \in F$ by \[\alpha(x,y) = \min\set{v(x) - \frac 1 2 v(b), v(y) + \frac 1 2 v(b)}\] is a $v$-norm which is 
compatible of depth $\varepsilon = 0$ with $q$. \\
Finally,  if $v(a) + v(b) > 0$ and $v(a), v(b) \geq 0$, then the map $\alpha \colon F \times F \to \Gamma \cup \set{\infty}$ defined for all $x,y \in F$ by \[\alpha(x,y) = \min\set{v(x), v(y)}\] is also a $v$-norm which is compatible of depth $\varepsilon = 0$ with $q$. \
\
(3) Now suppose $\charac F = 0$. Then, the map $\alpha \colon F \to \Gamma \cup \set{\infty}$ defined by \[\alpha(x) = \frac 1 2 v(a) + v(x)\] for all $x \in F$ is clearly a $v$-norm which is compatible of depth \(\varepsilon = v(2)\) with $q$. 
}

\prop{\label{existence_of_norms}Every nonsingular quadratic space admits a compatible $v$-norm of depth $\varepsilon$, for some \(\varepsilon \in \frac 1 2 \Gamma_F\) such that \(0 \leq \varepsilon \leq v(2)\). Moreover, every hyperbolic quadratic space 
admits a tame compatible \(v\)-norm.
}

\pr{Let $q \colon V \to F$ be a nonsingular quadratic form with polar form $b_{q}$. Assume first that $\charac F = 2$ (so that $v(2) = \infty$) and $\dim V = 2$. Since $b_q$ is nondegenerate, there exist $e, f \in V$ such that $b_q(e, f) = 1$, so that $q$ is 
isometric to the form $(x,y) \in F \times F \mapsto q(e)x^2+ xy + q(f)y^2$ which admits a compatible $v$-norm of depth $\varepsilon \geq 0$ by Lemma~\ref{existence_of_norms_lemma}. Now suppose $\charac F = 0$. In this case, we have $0 \leq v(2) < \infty$. Suppose also $\dim V = 1$ and pick $e \in V$ so that $q(e) \neq 0$. Then, $q$ is isometric to the form $x \in F \mapsto q(e)x^2$ which admits a compatible $v$-norm of depth $\varepsilon = v(2)$ by Lemma~\ref{existence_of_norms_lemma}. In 
general, no matter the characteristic of $F$, if $\dim V = n > 0$, we write $V = (V_1, q_1) \perp \dots \perp (V_n, q_n)$ for some $F$-quadratic spaces $(V_i, q_i)$  for $i = 1, \dots, n$, each of them admitting a compatible $v$-norm $\alpha_i$ of depth $\varepsilon_i \leq v(2)$. Note that such a decomposition exists by the first two steps and by normalisation of nonsingular quadratic forms, cf. Lemma~\ref{normalisation_gr} applied to \(b_q\) for \(F\) with trivial graduation. Then, by Lemma~\ref{depth_of_norms}, 
we can assume that all the $\alpha_i$'s are compatible of the same depth \(\varepsilon = \max_{i=1, \dots,n}\varepsilon_i \in \frac 1 2 \Gamma_F\) with \(0 \leq \varepsilon \leq v(2)\). Finally, by Lemma~\ref{sum_of_norms}, the map $\alpha = \alpha_1 \oplus \dots \oplus \alpha_n$ is a compatible $v$-norm of depth $\varepsilon$. For the second part of the proof, if $(V,q)$ is moreover hyperbolic, then it can be decomposed as an orthogonal sum of hyperbolic planes, each of them admitting a compatible $v$-norm $\alpha_i$ of depth $0$ by Lemma~\ref{existence_of_norms_lemma}. This concludes the proof.
}

\lem{\label{hyperMeta}Let $\varphi = (V,q)$ be a quadratic space with a compatible $v$-norm $\alpha$ of depth $\varepsilon$. If $\varphi$ is hyperbolic, then $\widetilde{\varphi}_{\alpha}$ is metabolic.
}

\pr{If $U \subset V$ is a totally isotropic subspace for $q$ of dimension $\frac 1 2 \dim V$, then $\gr_{\alpha}(U)$ is a Lagrangian for $\widetilde{\varphi}_{\alpha}$.
}

\lem{\label{depth_Witt}Let $\varphi = (V,q)$ be a quadratic space and $\alpha, \beta$ two $v$-norms which are compatible with $q$ of the same depth $\varepsilon$. Then the spaces $\widetilde{\varphi}_{\alpha}$ and $\widetilde{\varphi}_{\beta}$ are Witt 
equivalent.
}

\pr{By Lemma~\ref{sum_of_norms}, the space $(- \varphi) \perp \varphi$ admits $\alpha \oplus \beta$ as compatible $v$-norm of depth $\varepsilon$. Therefore, by Lemma~\ref{hyperMeta}, the graded space $(- \widetilde{\varphi}_{\alpha}) \perp 
\widetilde{\varphi}_{\beta}$ is metabolic.
}

\thm{\label{meta_depth}Let $\alpha$ be a $v$-norm which is compatible of depth $\gamma \in \Gamma$ with a quadratic space $\varphi = (V,q)$. If $0 < \gamma \leq v(2)$, the space $\varphi$ admits a compatible $v$-norm of depth $\delta < \gamma$ if and 
only if the space $\widetilde{\varphi}_{\alpha}$ is metabolic. 
}

\pr{First suppose that $\varphi$ admits a compatible norm $\beta'$ of depth $\delta < \gamma \leq v(2)$. Then by Lemma~\ref{depth_of_norms}, there exists a norm $\beta = \beta' - \frac 1 2 (\gamma - \delta)$ of depth $\gamma$ such that
\[
\begin{aligned}
	&\beta'(x) + \beta'(y) + \delta = \beta(x) + \beta(y) + \gamma &\text{ and } &\beta(x) < \beta'(x).
\end{aligned}
\]
We will see that $\widetilde{\varphi}_{\beta}$ is metabolic, hence $\widetilde{\varphi}_{\alpha}$ is metabolic by Lemma~\ref{depth_Witt}. Since $\delta < v(2)$, $\widetilde{(b_q)}_{\beta'}$ is an alternating bilinear form. Consequently, there exists a symplectic 
basis $\widetilde{e}_1, \widetilde{f}_1,\dots, \widetilde{e}_n, \widetilde{f}_n$ of $\widetilde{\varphi}_{\beta}$ such that $v(b_q(e_i,e_j)) > \beta'(e_i) + \beta'(e_j) + \delta$ for $i,j = 1, \dots, n$. Therefore $v(b_q(e_i,e_j)) > \beta(e_i) + \beta(e_j) + \gamma$ for all 
$i,j=1, \dots, n$. Since $e_1, f_1\dots, e_n, f_n$ is a splitting basis for $\beta'$, it is also a splitting basis of $\beta$ (by construction of $\beta$), hence $\widetilde{e}_1, \widetilde{f}_1, \dots, \widetilde{e}_n, \widetilde{f}_n \in \gr_{\beta}(V)$ is a basis such that 
the graded subspace spanned by $\widetilde{e}_1, \dots, \widetilde{e}_n$ is a totally isotropic subspace for $\widetilde{b}_{\beta}$. Moreover, the condition $\beta(x) < \beta'(x)$ implies that $\widetilde{q}_{\beta}$ is identically zero. That shows $
\widetilde{\varphi}_{\beta}$ is metabolic. Conversely, assume that $\widetilde{\varphi}_{\alpha}$ is metabolic and decompose it as a sum of metabolic planes, by Proposition~\ref{decomposition_gr}. We can therefore find a basis $e_1, f_1, \dots, e_n,f_n$ of 
$V$ which splits $\alpha$ and such that for all $i,j=1,\dots,n$ with $i\neq j$, 
\[
\begin{aligned}
	&v(b_q(e_i,e_j)) > \alpha(e_i) + \alpha(e_j) + \gamma \text{, } &v(b_q(e_i,f_i)) = \alpha(e_i) + \alpha(f_i) + \gamma, \\
  &v(b_q(e_i,f_j)) > \alpha(e_i) + \alpha(f_j) + \gamma \text{, } &v(b_q(f_i,f_j)) > \alpha(f_i) + \alpha(f_i) + \gamma, \\
  &v(q(e_i)) > 2\alpha(e_i) \text{, } &v(q(f_i))\geq 2 \alpha(f_i).
\end{aligned}
\]
We put 
\begin{multline*}
  \varepsilon = \min \bigg\{ \frac 1 2 v(q(e_i)) - \alpha(e_i), v(b_q(e_i,e_j))-\alpha(e_i)-\alpha(e_j) - \gamma, \\ \gamma \tq i,j=1, \dots, n \text{ and } i \neq j \bigg\} \in \Gamma,
\end{multline*}
($0 < \varepsilon \leq \gamma$), and we define a new $v$-norm $\alpha'$ on $V$ by
\[
  \alpha' \left( \sum_{i=1}^n \left( \lambda_i e_i + \mu_i f_i\right) \right) =
  \min \set{\alpha(\lambda_i e_i) + \varepsilon, \alpha(\mu_i f_i) \tq i=1, \dots, n}.
\]
It is easily seen that $\alpha'$ is compatible of depth $\gamma - \varepsilon$ with $q$.
}

For each $\varepsilon \in \Gamma$ such that $0 \leq \varepsilon \leq v(2)$, let $W_g(F, \varepsilon) : = W_{T}^{\varepsilon}(\gr_{v}(F))$ be the Witt group of $\varepsilon$-shifted spaces of type $T$ over \(\gr_v(F)\), where $T = I$ if $\varepsilon = 0$, $T = II$ if 
$0 < \varepsilon < v(2)$, and $T = \widetilde{2}-III$ if $\varepsilon = v(2)$. Therefore,  
\[
  W_g(F, \varepsilon) \simeq
  \begin{cases}
  W_q(\gr_{v}(F)) & \quad \text{if } \varepsilon = 0 \\
  W_{sq}^{\varepsilon}(\gr_{v}(F)) & \quad \text{if } 0 < \varepsilon < v(2) \\
  W(\gr_{v}(F)) & \quad \text{if } \varepsilon = v(2) \\
  \end{cases}
\]
Note that if $v(2) = 0$, $W_g(F, \varepsilon)$ is well-defined since in this case $0$-shifted spaces of type $I$ are exactly $v(2)$-shifted spaces type $\widetilde{2}-III$. Note also that \(W_g(F, \varepsilon) = 0\) when \(\varepsilon \not \in \frac 1 2 \Gamma_F\) by 
Lemma~\ref{origin_of_the_cases}. Let $\varepsilon \in \Gamma$ be such that $0 \leq \varepsilon \leq v(2)$ and let $W_{q}(F)_{\varepsilon} \subset W_{q}(F)$ be the set of Witt classes represented by a quadratic space which admits a compatible $v$-norm of 
depth $\gamma \leq \varepsilon$ (hence also a compatible $v$-norm of depth $\varepsilon$, by Lemma~\ref{depth_of_norms}). Since a hyperbolic space always admits a tame compatible norm (by Proposition~\ref{existence_of_norms}) and by 
Lemma~\ref{sum_of_norms}, $W_{q}(F)_{\varepsilon}$ is a subgroup of $W_{q}(F)$. Moreover, by definition, if $\gamma \leq \varepsilon$ then $W_{q}(F)_{\gamma} \subset W_{q}(F)_{\varepsilon}$. Therefore, by Proposition~\ref{existence_of_norms}, we get 
an ascending filtration by subgroups of $W_{q}(F)$ such that 
\[
   \begin{aligned}
      &\bigcup_{\varepsilon \in E} W_q(F)_{\varepsilon} = W_q(F), &\text{ where } E = \left\{ \varepsilon \in \frac 1 2 \Gamma_F \, | \, 0 \leq \varepsilon \leq v(2) \right\}.
   \end{aligned}
\]
Note that if $v(2) \in \Gamma$, then we have $W_q(F)_{v(2)} = W_q(F)$ since in that case \[W_q(F)_{v(2)} = \bigcup_{\varepsilon \in E} W_q(F)_{\varepsilon} = W_q(F).\]

Note also that for our valued field $F$, Proposition~\ref{normalisation_gr} induces a decomposition into two-dimensional spaces of every representative $\varphi$ of a class $[\varphi] \in W_g(F, \varepsilon)$ when $\varepsilon < v(2)$. But when $\varepsilon = 
v(2)$, the same proposition induces an orthogonal decomposition of every anisotropic representative $\varphi$ of a class $[\varphi] \in W_g(F, \varepsilon)$.

\thm{\label{filtration}Let $\varepsilon \in \Gamma$ such that $0 \leq \varepsilon \leq v(2)$. There exists a group epimorphism
\[
  \partial^{\varepsilon} \colon W_q(F)_{\varepsilon} \to W_g(F,\varepsilon)
\]
that carries the Witt class of a nonsingular quadratic space $\varphi$ with a compatible $v$-norm $\alpha$ of depth $\varepsilon$ to the Witt class of $\widetilde{\varphi}_{\alpha}$. If $\varepsilon > 0$, the kernel of this map is given by $\ker \partial^{\varepsilon} 
= W_{q}(F)_{< \varepsilon}$, which is the subgroup of $W_{q}(F)$ consisting of Witt classes with a representative admitting a compatible $v$-norm of depth $\gamma < \varepsilon$. If $\varepsilon = 0$ and $F$ is Henselian, then $\partial^{\varepsilon}$ is an 
isomorphism.
}

Note that this result was already proved by Elomary and Tignol for the subgroup \(W_q(F)_{0}\) in \cite{MR2836073}*{Proposition 8 and Theorem 10}.

\pr{Using Lemma~\ref{sum_of_norms}, Lemma~\ref{hyperMeta} and Lemma~\ref{depth_Witt}, it is routine to check that the map $\partial^{\varepsilon} \colon W_q(F)_{\varepsilon} \to  W_g(F,\varepsilon)$ is a well-defined group homomorphism. Suppose first \(\varepsilon < v(2)\). In this case, in order to prove the surjectivity of \(\partial^{\varepsilon}\), we only need to check that every anisotropic $2$-dimensional $\varepsilon$-shifted space $(V', q', b')$ is in the image. Here the proof uses the same ideas as 
\cite{MR2836073}*{Proposition 8}. Since $b'$ is nondegenerate and $q'$ anisotropic, there exist homogeneous $\xi_1, \xi_2 \in V'$ and $a_1,a_2 \in F^{\times}$ such that $b'(\xi_1,\xi_2) = 1$,  $q'(\xi_1) = \widetilde{a_1}$ and $q'(\xi_2) = \widetilde{a_2}$ (so 
that $\frac 1 2 v(a_1) + \frac 1 2 v(a_2) + \varepsilon = 0$). Consider the quadratic form given for all $x_1,x_2 \in F$ by $q(x_1, x_2) = a_1 x_1^2 + x_1 x_2 + a_2 x_2^2$, which is nondegenerate since $\varepsilon < v(2)$. Then $(\gr_{\alpha}(F \times F), 
\widetilde{q}_{\alpha}, \widetilde{b}_{\alpha}) \simeq (V', q', b')$ under the map $(x_1, x_2) \mapsto \widetilde{x_1}\xi_1 + \widetilde{x_2}\xi_2$, where the compatible $v$-norm $\alpha$ of depth $\varepsilon$ is given for all nonzero vector $(x_1, x_2) \in F 
\times F$ by $\alpha(x_1, x_2) = \min \set{\frac 1 2 v(a_1) + v(x_1), \frac 1 2 v(a_2) + v(x_2)}$. If $\varepsilon = v(2)$, in order to prove the surjectivity of \(\partial^{\varepsilon}\), we only need to check that every anisotropic one-dimensional $\varepsilon$-
shifted space $(V', q', b')$, is in the image. Since $b'$ is nondegenerate, we may find a homogeneous vector $\xi \in V$ such that $b'(\xi,\xi) \neq 0$. Consequently, there exists $a \in F^{\times}$ such that $q'(\xi) = \widetilde{a}$. Now consider the quadratic 
form on $F$ given for all $x\in F$ by $q(x) = a x^2$, which is nondegenerate when $v(2) \in \Gamma$, and define the $v$-norm $\alpha \colon F \to \Gamma \cup \set{\infty}$ given for every nonzero vector $x \in F$ by $\alpha(x) = \frac 1 2 v(a) + v(x)$, which 
is compatible with $(F, q)$ of depth $\varepsilon = v(2)$. From this, straightforward computations show that $(\gr_{\alpha}(F), \widetilde{q}_{\alpha}, \widetilde{b}_{\alpha}) \simeq (V', q', b')$ under the map $x \mapsto \widetilde{x}\xi$. When $\varepsilon > 0$, 
the identity $\ker \partial^{\varepsilon} = W_{q}(F)_{< \varepsilon}$ follows from Theorem~\ref{meta_depth}. When $\varepsilon = 0$ and $F$ is Henselian, the injectivity follows from \cite{MR2836073}*{Theorem 10}.
}

Note that the Henselian hypothesis is required only for \(\partial^{0} \colon W_q(F)_{0} \to W_g(F,0)\) being injective.  That \(\partial^{\varepsilon}\) induced an isomorphism \(W_q(F)_{\varepsilon}/ W_{q}(F)_{< \varepsilon} \iso W_g(F,\varepsilon)\) if \(\varepsilon 
> 0\) holds for every valued field \((F,v)\).

\subsection{\texorpdfstring{$W_{q}(F)_{\varepsilon}$}{} with generators}

In this section, we give another description of the subgroups $W_q(F)_{\varepsilon}$ when $0 \leq \varepsilon < v(2)$, and we complete the proof of the result $(2)$ of the introduction by showing that 
\[
   \begin{aligned}
      &W_{q}(F)_{< \varepsilon} = \bigcup_{\substack{\gamma \in E, \\ 0 \leq \gamma < \varepsilon}} W_q(F)_{\gamma}, &\text{ where } E = \set{\gamma \in \frac 1 2 \Gamma_F \tq 0 \leq \gamma \leq v(2)}
   \end{aligned}
\]
 if \(0 < \varepsilon \leq v(2)\). For that, we first prove the next proposition.

\prop{\label{normalisation_espaceNorme}If $(V,q)$ is a quadratic space and $\alpha$ is a $v$-norm on $V$ that is compatible of depth $\varepsilon < v(2)$ with $q$, then we can find subspaces $(V_1, q_1), \dots, (V_n, q_n)$ of dimension $2$ such that
\listeN{
  \item $(V,q) = (V_1,q_1) \perp \dots \perp (V_n, q_n)$, and
  \item each $\alpha_i = \alpha_{|V_i}$ is a $v$-norm on $V_i$ which is \(\varepsilon\)-compatible with $q_i$, and $\alpha = \alpha_1 \oplus \dots \oplus \alpha_n$.
}
}

For the proof, we will use the following result, which is a particular case of \cite{coyette}*{Proposition 2.5}.

\prop{\label{Coyette}Let $(F,v)$ be a valued field (of arbitrary characteristic). Consider $(V, \alpha)$ a $F$-vector space with $v$-value function $\alpha$, and a subspace $U \subset V$. Suppose that $p \colon V \to U$ is a linear map such that $p(u) = u$ for 
all $u \in U$. If $\alpha(p(x)) \geq \alpha(x)$ for all $x \in V$, then $V = U \oplus \ker p$ and $\alpha(u + v) = \min \set{\alpha(u), \alpha(v)}$ for all $u \in U$ and $v \in \ker p$. 
}
\begin{proof}[Proof of Proposition~\ref{normalisation_espaceNorme}] Let $(V, q)$ be a quadratic space admitting a compatible $v$-norm $\alpha$ of depth $\varepsilon$, with $0 \leq \varepsilon < v(2)$. Suppose $\dim V > 0$ and take two nonzero vectors $e,f 
\in V$ such that $v(b_{q}(e,f)) = \alpha(e) + \alpha(f) + \varepsilon$. Consider the subspace $U = \sev \langle e, f \rangle \subset V$ spanned by $e$ and $f$. Since $\varepsilon < v(2)$, the vectors $e$ and $f$ are linearly independent, and ${b_q}_{|U}$ is 
nondegenerate. So we can write \[ \left(V, q\right) = \left(U, q_{|U}\right) \perp \left(U^{\perp}, q_{|U^{\perp}}\right) \] for the two nonsingular quadratic spaces $(U, q_{|U})$ and $(U^{\perp}, q_{|U^{\perp}})$. Since $\alpha$ is a $v$-norm, it follows that $\alpha_{|
U}$ and $\alpha_{|U^{\perp}}$ are $v$-norms too (see \cite{MR3328410}*{Proposition 3.14}). Since \(\varepsilon < v(2)\), it is straightforward to see that \(\alpha_{|U}\) is \(\varepsilon\)-compatible with $(U, q_{|U})$. To conclude the proof, it remains to show that 
$\alpha_{|U^{\perp}}$ is \(\varepsilon\)-compatible with $(U^{\perp}, q_{|U^{\perp}})$, because the proposition then follows by induction on the dimension. Consider the orthogonal projection $p \colon V \to U$ given for all $x \in V$ by 
\[
	p(x) = \frac{\Delta_1(x)}{\Delta_0} e + \frac{\Delta_2(x)}{\Delta_0} f,
\]
where 
\[
   \begin{aligned}
      \Delta_0 &= \det \begin{pmatrix} b_q(e,e) & b_q(e,f) \\ b_q(e,f) & b_q(f,f) \end{pmatrix}, & \qquad
      \Delta_1(x) &=  \det \begin{pmatrix} b_q(x,e) & b_q(e,f) \\ b_q(x,f) & b_q(f,f)  \end{pmatrix}, 
   \end{aligned}
\]
\[
   \Delta_2(x) = \det \begin{pmatrix} b_q(e,e) & b_q(e,x) \\ b_q(e,f) & b_q(f,x) \end{pmatrix}.
\]
Using that 
\[ 
   v(b_{q}(e,f)) = \alpha(e) + \alpha(f) + \varepsilon \text{ and } \varepsilon < v(2),
\] 
standard calculations yield 
\[
   \begin{aligned}
      v\left(\Delta_0\right) = 2 \alpha(e) + 2 \alpha(f) + 2 \varepsilon, & \qquad
      v\left(\Delta_1(x)\right) \geq \alpha(e) + 2 \alpha(f) + \alpha(x) + 2 \varepsilon,
   \end{aligned}
\]
\[
   v\left(\Delta_2(x)\right) \geq 2 \alpha(e) + \alpha(f) + \alpha(x) + 2 \varepsilon.
\]
Those inequalities imply $\alpha(p(x)) \geq \alpha(x)$. Therefore, by Proposition~\ref{Coyette},  \[\alpha(u + v) = \min \set{\alpha(u), \alpha(v)}\] for all $u \in U$ and $v \in U^{\perp}$. Now, suppose $x \in U^{\perp}$ is a nonzero vector. There exists a nonzero $
(y_1, y_2) \in U \oplus U^{\perp}$ such that \[v(b_{q}(x, y_1 + y_2)) = \alpha(x) + \alpha(y_1 + y_2) + \varepsilon.\] This implies: 
\[
\begin{aligned}
 v(b_{q}(x,y_2)) &\geq \alpha(x) + \alpha(y_2) + \varepsilon \\
 &\geq \alpha(x) + \min \set {\alpha(y_1), \alpha(y_2)} + \varepsilon = v(b_{q}(x, y_1 + y_2)) = v(b_{q}(x,y_2)).
\end{aligned}
\]
Consequently $v(b_{q}(x,y_2)) = \alpha(x) + \alpha(y_2) + \varepsilon$ for some nonzero $y_2 \in U^{\perp}$, and $\alpha_{|U^{\perp}}$ is  compatible with $(U^{\perp}, q_{|U^{\perp}})$ of depth $\varepsilon$.
\end{proof}

For $a, b \in F$, we denote by $[ a, b ]$ the quadratic space $(F \times F, q)$ where $q$ is given by $q(x_1,x_2) = a x_1^2 + x_1 x_2 + b x_2^2$ for all $x_1,x_2 \in F$. The following description of $W_q(F)_{\varepsilon}$ is inspired by the definitions of 
\cite{RH-05-2016}*{section 2}.

\cor{\label{W_q__asGenerated}Let $\varepsilon \in \Gamma$ such that $0 \leq \varepsilon < v(2)$. Then $W_q(F)_{\varepsilon}$ is the subgroup generated by the classes represented by some form $[ a, b ]$ with $a, b \in F$ and $v(a) + v(b) \geq - 
2\varepsilon$.
}

\pr{First note that, for $\varepsilon < v(2)$, the condition $v(a) + v(b) \geq - 2\varepsilon$ implies that $[a, b]$ is nonsingular. Suppose now $[ a, b ]$ is a quadratic space with $a, b \in F$ and $v(a) + v(b) \geq - 2\varepsilon$. Then, by 
Lemma~\ref{existence_of_norms_lemma}, $[a, b]$ admits a compatible $v$-norm of depth $\leq \varepsilon$. Conversely, assume $q \colon V \to F$ is a $2$-dimensional nonsingular quadratic form that admits a compatible $v$-norm of depth $\gamma \leq 
\varepsilon$. Pick $e,f \in V$ such that $v(b_q(e,f)) = \alpha(e) + \alpha(f) + \gamma$ and moreover $b_q(e,f) = 1$. Then \(e,f\) are linearly independent since \(\gamma < v(2)\), and $(V,q) \simeq [ q(e), q(f) ]$ with \[v(q(e)) + v(q(f)) \geq - 2 \gamma \geq  - 2 
\varepsilon.\] We have thus shown that the spaces $[a, b]$ with $a, b \in F$ and $v(a)+ v(b) \geq -2\varepsilon$ are the $2$-dimensional nonsingular spaces that carry a compatible $v$-norm of depth $\leq \varepsilon$, and the corollary follows by 
Proposition~\ref{normalisation_espaceNorme}. 
}

In order to prove completely the result announced in (2) of the introduction, it remains to prove the following.

\lem{Let $\varphi=(V,q)$ be a quadratic space. If $\varphi$ admits a compatible $v$-norm of depth $\varepsilon$ such that $0 \leq \varepsilon \leq v(2)$, then $\varphi$ also admits a compatible $v$-norm of depth $\gamma \leq \varepsilon$ such that $0 \leq 
\gamma \leq v(2)$ and moreover $\gamma \in \frac 1 2 \Gamma_F$. Consequently, 
\[
   \begin{aligned}
      &W_{q}(F)_{< \varepsilon} = \bigcup_{\substack{\gamma \in E, \\ 0 \leq \gamma < \varepsilon}} W_q(F)_{\gamma}, &\text{ where } E = \set{\gamma \in \frac 1 2 \Gamma_F \tq 0 \leq \gamma \leq v(2)}.
   \end{aligned}
\]
}

\pr{We only show the existence of the $v$-norm of depth $\gamma \leq \varepsilon$ such that $\gamma \in \frac 1 2 \Gamma_F$. If $\varepsilon = v(2)$, the assertion is clear. Suppose now $0 \leq \varepsilon < v(2)$. First assume $\varphi$ is a $2$-
dimensional space. Then, by the proof of Corollary~\ref{W_q__asGenerated}, write $\varphi \simeq [a,b]$ for some $a,b \in F$ such that $v(a) + v(b) \geq -2\varepsilon$. Therefore, by Lemma~\ref{existence_of_norms_lemma}, $[a, b]$ admits a compatible 
$v$-norm of depth $ \gamma \leq \varepsilon$ with $\gamma \in \frac 1 2 \Gamma_F$. Now, if $\varphi$ is a general nonsingular quadratic space which admits a compatible $v$-norm of depth $\varepsilon < v(2)$, decompose it, by 
Proposition~\ref{normalisation_espaceNorme}, into $2$-dimensional subspaces $\varphi_i$, each of them admitting a compatible $v$-norm of depth $\varepsilon$. By the first part of the proof, each space $\varphi_i$ admits a compatible $v$-norm of depth $
\gamma_i \leq \varepsilon$ for some $\gamma_i \in \frac 1 2 \Gamma_F$. Therefore, by Lemmas~\ref{depth_of_norms} and \ref{sum_of_norms}, $\varphi$ admits a compatible $v$-norm of depth $\gamma = \max_i \gamma_i \in \frac 1 2 \Gamma_F$.
}

\section{Relation with Arason's results}
\label{section-5}

In this section, we relate our work in the particular case of $\Gamma_F = \Z$ with Arason's results, and we give an example of application of those results. From now on, we suppose $\Gamma_F = \Z$. Note that in this case, for $\varepsilon \in \frac 1 2 \Z$ 
such that $0 < \varepsilon \leq v(2)$, we have $W_{q}(F)_{< \varepsilon} = W_q(F)_{\varepsilon - (1/2)}$. 

Let $\varphi$ be a nonsingular quadratic space. Since \(\Gamma_F = \Z\), there exists a minimal depth $\varepsilon \geq 0$ such that $\varphi$ admits a compatible $v$-norm of depth $\varepsilon$. We put $w(\varphi) = \varepsilon$ and we call it the 
\emph{wildness index of $\varphi$}. Let $\pi \in F$ be such that $v(\pi) = 1$. Note that the ascending filtration $(W_q(F)_{\varepsilon})_{\varepsilon \in E}$, where $E = \set{\varepsilon \in \frac 1 2 \N \tq 0 \leq \varepsilon \leq v(2)}$, is infinite if $\charact F = 2$, 
whereas it is finite (and satisfies $W_q(F)_{v(2)} = W_q(F)$) if $\charact F \neq 2$.

Observe that, in Arason's note~\cite{RH-05-2016}, the index \(\varepsilon\) of the subgroups \(W_q(F)_{\varepsilon}\) ranges from \(0\) to \(2v(2)\). In this paper, we use half the values, so that \(\varepsilon\) ranges from \(0\) to \(v(2)\).

\subsection{Arason's isomorphisms}

In his note~\cite{RH-05-2016}, Arason defines the subgroups $W_{q}(F)_{\varepsilon}$, for a discretely valued field $F$, in terms of $2$-dimensional generators. He also describes those subgroups using isomorphisms which are given in terms of those 
generators. Corollary~\ref{Arason_2} is the corresponding description explained from our point of view.

\cor[compare \cite{RH-05-2016}*{section 2}]{\label{Arason_1}For $\varepsilon \in \frac 1 2 \N$ such that  $\varepsilon < v(2)$ and $n = 2 \varepsilon$, we have that $W_q(F)_{\varepsilon}$ is the subgroup generated by the classes represented by some form $
[ \alpha, \pi^{-n} \beta ]$ or $ \pi [ \alpha, \pi^{-n} \beta ] \simeq [ \pi \alpha, \pi^{-1-n} \beta ]$, with $\alpha, \beta \in F$ such that $v(\alpha) \geq 0$ and $v(\beta) \geq 0$.}

\pr{By Corollary~\ref{W_q__asGenerated}, $W_q(F)_{\varepsilon}$ is the subgroup generated by the classes represented by some forms $[a,b]$ with $a,b \in F$ such that $v(a) + v(b) \geq -2\varepsilon$. Note that, given an arbitrary $\pi \in F^{\times}$, we 
have $\pi [ a, b ] \simeq [ \pi a, \pi^{-1} b ]$ and $\pi^2 [ a, b ] \simeq [ a, b ]$ for all $a, b \in F$. To complete the proof, it remains to apply these general relations with \(\pi \in F^{\times}\) such that \(v(\pi) = 1\). 
}

For the rest of the paper, we write $\varphi_W$ for the Witt class in $W_q(F)$ of a nonsingular quadratic space $\varphi$. Given elements \(\alpha_i \in F\) (\(i = 1, \dots, n\)), we denote by $\langle \alpha_1, \dots, \alpha_n \rangle$ the symmetric bilinear space 
$(V,b)$ with orthogonal basis \(e_1, \dots, e_n \in V\) such that $b(e_i,e_i) = \alpha_i$ and $b(e_i,e_j) = 0$ if $i \neq j$, for all \(i,j = 1, \dots, n\). When $\charact F \neq 2$, we identify nondegenerate symmetric bilinear forms with nonsingular quadratic forms 
(with a symmetric bilinear form $b$, we associate the quadratic form given by $x \mapsto b(x,x)$, and with a quadratic form $q$ we associate $\frac 1 2 b_q$).

\cor[\cite{RH-05-2016}*{Proposition 1.1, Proposition 3.1 and Proposition 2.1}]{\label{Arason_2} Assume\footnote{Observe that if $v(2) = 0$, the well-known epimorphism (or isomorphism when $F$ is Henselian) $W_q(F) \left(=W_q(F)_{v(2)}= W_q(F)_0\right) 
\to W_q(\overline{F}) \oplus W_q(\overline{F})$ can also be deduced from our results (and is also in \cite{RH-05-2016}*{Proposition 1.1}).} $\charact \overline{F}= 2$. In the following, $\alpha, \beta$ represent elements in $F$ such that $v(\alpha) \geq 0$ and 
$v(\beta) \geq 0$.
There is always a group epimorphism
\[
  W_q(F)_0 \to W_q(\overline{F})  \oplus W_q(\overline{F}) 
\]
which maps a Witt class $[ \alpha,  \beta ]_W$ to $([\overline{\alpha}, \overline{\beta}]_W,0)$ and a Witt class $[\pi \alpha, \pi^{-1}\beta]_W$ to $(0,[\overline{\alpha}, \overline{\beta}]_W)$. That epimorphism is an isomorphism when $F$ is Henselian. \\
Let $\varepsilon \in \frac 1 2 \N$ such that $0 < \varepsilon < v(2)$. 
If $\varepsilon \notin \N$, there exists a group isomorphism 
\[
  W_q(F)_{\varepsilon} / W_q(F)_{\varepsilon-(1/2)} \to \overline{F} \otimes_{\overline{F}^2} \overline{F} 
\]
which maps the class of a Witt class $[ \alpha, \pi^{-2\varepsilon} \beta ]_W$ to $\overline{\alpha} \otimes \overline{\beta}$ and the class of a Witt class $[\pi \alpha, \pi^{-1-2\varepsilon} \beta ]_W$ to $\overline{\beta} \otimes \overline{\alpha}$. \\
If $\varepsilon \in \N$, there exists a group isomorphism 
\[
  W_q(F)_{\varepsilon} / W_q(F)_{\varepsilon-(1/2)} \to \left( \overline{F} \wedge_{\overline{F}^2} \overline{F} \right) \oplus \left( \overline{F} \wedge_{\overline{F}^2} \overline{F} \right)
\]
which maps the class of a Witt class $[\alpha, \pi^{-2\varepsilon} \beta ]_W$ to $(\overline{\alpha} \wedge \overline{\beta}, 0)$ and the class of a Witt class $[\pi \alpha, \pi^{-1-2\varepsilon} \beta ]_W$ to $(0, \overline{\alpha} \wedge \overline{\beta})$. \\
If $\charact F \neq 2$, there is a group isomorphism 
\[
  W_q(F) / W_q(F)_{v(2) - (1/2)}  \to W(\overline{F}) \oplus W(\overline{F}), 
\] 
which maps the class of a Witt class $\langle \alpha \rangle_W$ with $v(\alpha) = 0$ to $(\langle \overline{\alpha} \rangle_{W}, 0)$ and the class of a Witt class $\langle \pi \beta \rangle_{W}$ with $v(\beta) = 0$ to $(0, \langle \overline{\beta} \rangle_{W})$.
}

\pr{Fix $\pi' \colon \frac 1 2 \Z / \Z \to \gr_{v}(F)$ such that $\pi'([0]) = 1$ and $\pi'([\frac 1 2]) = \widetilde{\pi}$. For the first map, compose the group homomorphisms of Theorem~\ref{filtration} and Proposition~\ref{decompsitionCan_1}, with a choice of 
uniformizing parameters given by $\mathcal C = (\rho, \pi')$ where $\rho = 1$. If $\varepsilon \notin \N$, compose the group homomorphisms of Theorem~\ref{filtration}, Proposition~\ref{decompsitionCan_1} and Theorem~\ref{structure_ssq}, with a choice of 
uniformizing parameters given by $\mathcal C = (\rho, \pi')$ where $\rho = \widetilde{\pi}^{2\varepsilon}$. If $\varepsilon \in \N$, compose the group homomorphisms of Theorem~\ref{filtration}, Proposition~\ref{iso_gr_sq_2} and Theorem~\ref{structure_sq}, with 
a choice of uniformizing parameters given by $\mathcal C = (\rho, \pi')$ where $\rho = \widetilde{\pi}^{\varepsilon}$. In those three cases (put $\varepsilon=0$ in the first), the map defined by $\alpha(x,y) = \min \set{v(x), v(y)-\varepsilon}$ for all $x,y \in F$ is a \(\varepsilon\)-compatible $v$-norm with $[ \alpha, \pi^{-2\varepsilon} \beta ]$, and the map defined by $\alpha(x,y) = \min \set{v(x) + \frac 1 2, v(y)-\varepsilon - \frac 1 2}$ for all $x,y \in F$ is a \(\varepsilon\)-compatible $v$-norm with $[\pi \alpha, 
\pi^{-1-2\varepsilon} \beta ]$. In the last case when $\charact F \neq 2$, compose the group homomorphism of Theorem~\ref{filtration} and Proposition~\ref{decompsitionCan_1}, with a choice of uniformizing parameters given by $\mathcal C = (\rho, \pi')$ 
where $\rho = \widetilde{2}$. Here the $v$-norms constructed in Lemma~\ref{existence_of_norms} are compatible of depth $v(2)$ with $\langle \alpha \rangle$ and $\langle \pi \beta \rangle$.
}

Corollary~\ref{Arason_2} can also be found in \cite{RH-19-2006}*{Theorem 2} in the particular case of $F = K((S))$ being the field of formal Laurent series in an indeterminate $S$ over a field $K$ of characteristic $2$.

\subsection{An application}

Our goal for this section is to prove Proposition~\ref{complete_discrete_dyadic} as an illustration of the results of Corollary~\ref{Arason_2}. In order to achieve this purpose, we need the following lemma.

\lem{\label{calculations_discrVal}Let $\varepsilon \in \frac 1 2 \N$ and $k \in \N$. For $a,b,c,d \in F$ and $\alpha, \beta, \gamma \in F$ such that $v(\alpha) \geq 0$, $v(\beta) \geq 0$ and $v(\gamma) \geq 0$, the following relations hold:
\begin{itemize}
\item[\((a)\)] \([a,b] \simeq [b,a]\)
\item[$(b)$]
   $[c^2 a, b] \simeq [a, c^2 b]$,
\item[$(c)$]
   If \(c \neq 0\) then \(c[a,b] \simeq [c^{-1}a, cb] \), hence 
   $\pi [\alpha, \beta \pi^{-2k-1}] \simeq [\beta, \alpha \pi^{-2k-1}]$.
\item[$(d)$] \emph{(\cite{RH-05-2016}*{Proposition 2.3})}
   If \([a,b]\) and \([c,d]\) are nonsingular then \[\left[a,b\right] \perp \left[c,d\right] \simeq \left[a + c, b\right] \perp \left[\frac{-c}{1-4cd}, \frac{d-b}{1-4ab}\right].\] 
   Hence if \(0 < \varepsilon < v(2)\) then
\[
\begin{aligned}
&[\alpha, \pi^{-2\varepsilon}\gamma]_W + [\beta, \pi^{-2\varepsilon}\gamma]_W \equiv [\alpha + \beta, \pi^{-2\varepsilon}\gamma]_W  \mod W_q(F)_{\varepsilon-(1/2)}, \\
 &[\alpha, \pi^{-2\varepsilon}\beta]_W + [\alpha, \pi^{-2\varepsilon}\gamma]_W \equiv [\alpha, \pi^{-2\varepsilon}(\beta+\gamma)]_W \mod W_q(F)_{\varepsilon-(1/2)}, \\
 &[\pi \alpha, \pi^{-1-2\varepsilon}\gamma]_W + [\pi \beta, \pi^{-1-2\varepsilon}\gamma]_W \equiv  [\pi(\alpha + \beta), \pi^{-1-2\varepsilon}\gamma]_W  \mod W_q(F)_{\varepsilon-(1/2)}, \\
&[\pi \alpha, \pi^{-1-2\varepsilon}\beta]_W + [\pi\alpha, \pi^{-1-2\varepsilon}\gamma]_W \equiv [\pi \alpha, \pi^{-1-2\varepsilon}(\beta+\gamma)]_W  \mod W_q(F)_{\varepsilon-(1/2)}. 
\end{aligned}
\]
 Moreover, if either \(\charact F = 2\) or \(F\) is a complete field satisfying \(\charact \overline{F} = 2\), then
 \[
  \begin{aligned}
  &[\alpha, \gamma]_W + [\beta, \gamma]_W = [\alpha + \beta, \gamma]_W, \\
  &[\alpha, \beta+\gamma]_W = [\alpha, \beta]_W + [\alpha, \gamma]_W, \\
  &[\pi \alpha, \pi^{-1}\gamma]_W + [\pi \beta, \pi^{-1}\gamma]_W = [\pi(\alpha + \beta), \pi^{-1}\gamma]_W, \\
  &[\pi \alpha, \pi^{-1}\beta]_W + [\pi\alpha, \pi^{-1}\gamma]_W = [\pi \alpha, \pi^{-1}(\beta+\gamma)]_W.
  \end{aligned}
\]
\item[$(e)$]
   $[a,b]_W = 0$ if $v(ab) > 0$  and $F$ is complete.
\item[$(f)$]
   If $\charact F \neq 2$ and $a \neq 0$, $\langle a, b \rangle \simeq [a, \frac{1}{4a^2}(a+b)]$
\item[$(g)$]
  If $\charact F \neq 2$, $\langle a \rangle_W + \langle b \rangle_W = \langle a + b \rangle_W + \langle ab(a+b) \rangle_W$ for all $a,b \in F^{\times}$ such that $a +b \neq 0$.
\end{itemize}
}

\pr{The easy proofs of $(a)$ and \((b)\) are left to the reader. The first part of $(c)$ is left to the reader, and the second part of \((c)\) follows directly from $(a)$ and \((b)\) since $\pi [\alpha, \pi^{-2k-1}\beta] \simeq [\pi^{-1}\alpha, \pi^{-2k}\beta]$ by the first part of 
\((c)\). For $(e)$, observe that since \(F\) is complete and \(v(ab) > 0\), by Hensel's lemma, there exists \(u \in F\) such that \(u^2 + u + ab = 0\). Now, if \(a = 0\), then $[a,b]_W = 0$. If \(a \neq 0\), \((a^{-1}u,1)\) is an isotropic vector for \([a,b]\) and the conclusion 
follows. The first part of $(d)$ comes from \cite{RH-05-2016}*{Proposition 2.3}. It is a direct calculation: if \(e_1,f_1,e_2,f_2\) is the standard basis of \(\left[a,b\right] \perp \left[c,d\right]\), then consider the change of basis given by \(e_1' = e_1 + e_2\), \(f_1' = 
f_1\), \(e_2' = (e_2 - 2cf_2)/(1-4cd)\) and \(f_2' = (2be_1 - f_1 +(1-4ab)f_2)/(1-4ab)\). The first part of \((d)\) yields \[\left[\alpha, \gamma \right] \perp \left[\beta, \gamma \right] \simeq \left[\alpha + \beta, \gamma \right] \perp \left[\frac{-\beta}{1-4\beta \gamma}, 
\frac{-2\gamma}{1-4\alpha \gamma}\right].\] The last summand is clearly isotropic if \(\charact F = 2\), hence \([\alpha, \gamma]_W + [\beta, \gamma]_W = [\alpha + \beta, \gamma]_W\) in this case. This last equality also holds if \(F\) is a complete field such that 
\(\charact \overline{F} = 2\) by \((e)\), since \(v(2) > 0\) in that case. Suppose now \(0 < \varepsilon < v(2)\), then \[\left[\frac{-\beta}{1-4\beta \pi^{-2\varepsilon}\gamma}, \frac{-2 \pi^{-2\varepsilon} \gamma}{1-4\alpha \pi^{-2\varepsilon} \gamma}\right] \in 
W_q(F)_{\varepsilon-(1/2)},\] by Corollary~\ref{W_q__asGenerated} since in this case we have \(v(4\beta\pi^{-2\varepsilon} \gamma) > 0\) and \(v(\alpha \pi^{-2\varepsilon} \gamma) > 0\), \(v(1-4\beta \pi^{-2\varepsilon}\gamma) = 0 \) and \(v(1-4\alpha 
\pi^{-2\varepsilon} \gamma) = 0\),  \(v(2) \geq 1\) and finally \(v(\beta 2\pi^{-2\varepsilon}\gamma) \geq -2(\varepsilon - \frac{1}{2})\). The rest of the proof of \((d)\) is left to the reader. For $(f)$, note that $\langle a, b \rangle \simeq [a, \frac{1}{4a^2}(a+b)]$ by the 
change of basis given by $e' = e$ and $f' = \frac{e-f}{2a}$, where $e,f$ is the standard basis of $\langle a, b \rangle$. Finally, $(g)$ is well-known. See for example \cite{MR2427530}*{Lemma 4.1}.
}

The following proposition and its proof are inspired by \cite{MR3437769}*{Lemma 8.1 and Proposition 8.2}. 

\prop{\label{complete_discrete_dyadic}Suppose $F$ is a complete discretely valued field with perfect dyadic residue field $\overline{F}$ and let $\pi \in F$ be such that $v(\pi) = 1$. If $\charact F = 2$, then every Witt class of a quadratic form over $F$ can be 
written as
\[
	\sum_{k=0}^{n} \left[1, \alpha_{2k+1}^2 \pi^{-2k-1}\right]_W + \left[1, \alpha_0^2\right]_W + \left[\pi, \beta_0^2\pi^{-1}\right]_W,
\]
for some $n \in \N$ and  $\alpha_k, \beta_k \in F$ such that $v(\alpha_k) \geq 0$ and $v(\beta_k) \geq 0$ for all $k$. If $\charact F \neq 2$, then every Witt class of a quadratic form over $F$ can be written as
\[
	\sum_{k=0}^{v(2)-1} \left[1, \alpha_{2k+1}^2 \pi^{-2k-1}\right]_W + \left[1, \alpha_0^2\right]_W + \left[\pi, \beta_0^2\pi^{-1}\right]_W + \alpha \langle 1 \rangle_W + \beta \langle \pi \rangle_W, 
\]
for some $\alpha, \beta \in \set{0,1}$ uniquely determined and some $\alpha_k, \beta_k \in F$ such that $v(\alpha_k) \geq 0$ and $v(\beta_k) \geq 0$ for all $k$. Moreover, in both cases, if we fix a (set) section $s \colon \overline{F} \to \mathcal O$ of the 
canonical quotient map $\mathcal O \to \overline{F}$ (where $\mathcal O := \set{x \in F \tq v(x) \geq 0}$), we can always choose $\alpha_k, \beta_k \in s(\overline{F})$ for all $k$. Provided that $\alpha_k, \beta_k \in s(\overline{F})$ for all $k$, the decomposition 
of a Witt class is unique.\footnote{When $\charact F = 2$, if $K$ is a coefficient field (i.e., a subfield $K$ of $F$ contained in $\mathcal{O}$ that maps isomorphically onto $\overline{F}$ under the canonical quotient map $\mathcal{O} \to \overline{F}$), one can 
choose for $s$ the isomorphism from $\overline{F}$ to $K$.}
}

Note that the decompositions given in Proposition~\ref{complete_discrete_dyadic} are in general neither decompositions with anisotropic forms nor the decompositions with the least number of summands. For example, in the case when \(\charact F = 2\), the 
stated decomposition can be rewritten as 
\[
    \left[1,  \sum_{k=0}^{n} \alpha_{2k+1}^2 \pi^{-2k-1} + \alpha_0^2\right]_W + \left[\pi, \beta_0^2\pi^{-1}\right]_W.
\]

\pr{First, let $\alpha,\beta \in F$ such that $v(\alpha) \geq 0$ and $v(\beta) \geq 0$. Since $\overline{F}^{2} = \overline{F}$, there exist $\alpha', \tau \in F$ such that $\alpha = \alpha'^2 + \tau$ and $v(\alpha') \geq 0$, $v(\tau) > 0$. Similarly, write $\beta = 
\beta'^2 + \mu$ for some $\beta', \mu \in F$ such that $v(\beta') \geq 0$ and $v(\mu) > 0$. Then, for $k \in \N$ such that $k < 2v(2)$, by Lemma~\ref{calculations_discrVal} $(b)$ and $(d)$, we have $[\alpha, \beta \pi^{-k}]_W = [1, (\alpha' \beta')^2 \pi^{-k}]_W + 
[ \alpha'^2 , \mu \pi^{-k}]  + [ \tau , \beta'^2 \pi^{-k}]  + [ \tau , \mu \pi^{-k}] + w$ for some $w \in W_q(F)$ such that $w \in W_q(F)_{(k-1)/2}$ if $k > 0$ and $w=0$ if $k = 0$. Observe that $[ \alpha'^2 , \mu \pi^{-k}], [ \tau , \beta'^2 \pi^{-k}], [\tau , \mu \pi^{-k}] \in 
W_q(F)_{(k-1)/2}$ if $k > 0$. Moreover, since $F$ is complete and $v(\alpha'^2 \mu), v(\tau \beta'^2), v( \tau \mu) > 0$, we have that $[ \alpha'^2 , \mu \pi^{-k}]$, $[ \tau , \beta'^2 \pi^{-k}]$ and $[ \tau , \mu \pi^{-k}]$ are hyperbolic when $k = 0$, by 
Lemma~\ref{calculations_discrVal} $(e)$. Similarly, $[\pi \alpha, \pi^{-1-k}\beta]_W = [\pi, \pi^{-1-k}(\alpha' \beta')^2] + w$ for some $w \in W_q(F)$ such that $w \in W_q(F)_{(k-1)/2}$ if $k > 0$ and $w=0$ if $k = 0$. Second, note that since $\overline{F}^{2} = 
\overline{F}$, we have $W_q(F)_{\varepsilon} \subset W_q(F)_{\varepsilon-(1/2)}$ when $\varepsilon \in \N\setminus \set{0}$. Indeed, that follows from Corollary~\ref{Arason_2} since here $\overline{F} \wedge_{\overline{F}^2} \overline{F} = \set{0}$. Third, we 
show the existence of the decomposition of a Witt class $\varphi_W$ by induction on the depth $\varepsilon$ of a $v$-norm which is \(\varepsilon\)-compatible with the form $\varphi$. Suppose that $\varphi$ admits a compatible $v$-norm of depth $
\varepsilon=0$. Then Corollary~\ref{Arason_1} and the first part of the proof show that $\varphi_W$ can be written as a sum of classes of two-dimensional spaces of the form $[1, a^2]$ and $[\pi, \pi^{-1}b^2]$ with $a,b \in F$ such that $v(a),v(b)\geq 0$. So we 
conclude by Corollary~\ref{calculations_discrVal} $(d)$ and the first part of the proof again. Suppose now that $\varphi$ admits a compatible $v$-norm of depth $\varepsilon \in \frac 1 2 \N$ satisfying $0 < \varepsilon < v(2)$. By the second part of the proof, we 
way assume $\varepsilon \notin \N$. By Corollary~\ref{Arason_1}, Lemma~\ref{calculations_discrVal} $(c)$ and the first part of the proof, $\varphi_W$ can be written as a sum of some $w \in W_q(F)_{\varepsilon - (1/2)}$ (hence $w \in W_q(F)_{\varepsilon - 
1}$, by the second part of the proof) and a sum of classes of two-dimensional spaces of the form $[1, a^2 \pi^{-2\varepsilon}]$ for some $a \in F$ such that $v(a)\geq 0$. So we conclude here by using by Corollary~\ref{calculations_discrVal} $(d)$, the first part 
of the proof, and the induction hypothesis. This shows the existence of the decomposition when $\charact F = 2$. Suppose now that $\charact F \neq 2$ and that $\varphi$ admits a compatible $v$-norm of depth $\varepsilon=v(2)$. Let $a \in F^{\times}$ be 
such that $v(a) = 0$ and write $a =\alpha^2 + \mu$ for some $\alpha \in F^{\times}$ such that $v(\alpha) = 0$ and $\mu \in F$ such that $v(\mu) > 0$.  If $\mu = 0$, then $\langle a \rangle \simeq \langle \alpha^2 \rangle \simeq \langle 1 \rangle$.  Otherwise, if $
\mu\neq 0$, then $\langle a \rangle_W = \langle \alpha^2 \rangle_W + \langle \mu \rangle_W - \langle \alpha^2 \mu a \rangle_W = \langle 1\rangle_W + \langle \mu, - \mu\alpha^{-2}a \rangle_W$, by Corollary~\ref{calculations_discrVal} $(g)$. But $\langle \mu, - 
\mu\alpha^{-2}a \rangle_W \in W_q(F)_{v(2)-(1/2)}$, since by Lemma~\ref{calculations_discrVal} $(f)$ $\langle \mu, - \mu\alpha^{-2}a \rangle  \simeq [\mu, 2^{-2}\mu^{-1}(1-\alpha^{-2}a)]$ with $v(\mu 2^{-2} \mu^{-1} (1-\alpha^{-2}a)) = - 2v(2) + v(1-\alpha^{-2}a) 
\geq -2(v(2)- \frac 1 2)$. Indeed, $v(1-\alpha^{-2}a))= v(\alpha^{-2} \mu) = v(\mu) \geq 1$. Note also that the space $\langle 1 ,  1 \rangle \simeq [2^{-1}, 1]$ admits a compatible norm of depth $\frac 1 2 v(2) < v(2)$. Since every nondegenerate bilinear form can 
be written as $\langle \alpha_1, \dots, \alpha_n \rangle \perp \langle \pi \rangle \langle \alpha'_1, \dots, \alpha'_m \rangle$ for some $\alpha_i, \alpha'_j$ such that $v(\alpha_i) = 0 = v(\alpha_j)$ for all $i,j$,  the existence of the decomposition in the second case 
can easily be concluded by the induction hypothesis. Now, the proof given a section $s \colon \overline{F} \to \mathcal O$ of the canonical quotient $\mathcal O \to \overline{F}$ is the same, except that each time we use the first part of the proof, we add the 
following step. We choose $\alpha'' \in s(\overline{F})$ such that $(\alpha' \beta')^2 = \alpha''^2 + \mu'$ for some $\mu \in F$ such that $v(\mu) > 0$. Then we write $[1, \pi^{-k} (\alpha' \beta')^2 ]_W = [1, \pi^{-k}\alpha''^2]_W + w$ for some $w \in W_q(F)$ such 
that $w \in W_q(F)_{(k-1)/2}$ if $k > 0$ and $w=0$ if $k = 0$. Similarly for $[\pi, \pi^{-1-k}(\alpha'\beta')^2]_W$. For the uniqueness, proceed as \cite{MR3437769}*{Proposition 8.2}. Note that in our case $W(\overline{F}) = \Z / 2 \Z$.
}

\subsection*{Acknowledgements}
This paper has been written during my doctoral training. I would like to thank my supervisor, Professor Jean-Pierre Tignol, for suggesting me this topic and for guiding me all along the road. Joachim Verstraete is a Research Fellow of the Fonds de la 
Recherche Scientifique - FNRS.

%%%%%%%%%%%%%%%%%%%%%%%%%%%%%%%%%%%%%%%%%%%%%%%%%%%%%%%%%%%%%%%%%%%%%%%%%%%%%%%%
%%%%%%%%%%%%%%%%%%%%%%%%%%%%%%%%%%%%%%%%%%%%%%%%%%%%%%%%%%%%%%%%%%%%%%%%%%%%%%%% 
%%%%%%%%%%%%%%%%%%%%%%%%%%%%%%%%%%%%%%%%%%%%%%%%%%%%%%%%%%%%%%%%%%%%%%%%%%%%%%%%

{\small 
\textsc{ICTEAM Institute/INMA, Université catholique de Louvain, Euler Building, Avenue Georges Lemaître 4 box L4.05.01, 1348 Ottignies-Louvain-la-Neuve, Belgium}\\
\indent \textit{E-mail address:} \verb+joachim.verstraete@uclouvain.be+
}

\end{document}